\documentclass[11pt]{amsart}
\usepackage{amsxtra}
\usepackage{amssymb,MnSymbol}
\usepackage{tikz-cd}
\theoremstyle{plain}

\newcommand{\cleqn}{\setcounter{equation}{0}}
\newcommand{\clth}{\setcounter{theorem}{0}}
\newcommand {\sectionnew}[1]{\section{#1}\cleqn\clth}
\newcommand{\nn}{\hfill\nonumber}
\newtheorem{theorem}{Theorem}[section]
\newtheorem{lemma}[theorem]{Lemma}
\newtheorem{definition-theorem}[theorem]{Definition-Theorem}
\newtheorem{proposition}[theorem]{Proposition}
\newtheorem{corollary}[theorem]{Corollary}
\newtheorem{definition}[theorem]{Definition}
\newtheorem{example}[theorem]{Example}

\newtheorem{remark}[theorem]{Remark}
\newtheorem{conjecture}[theorem]{Conjecture}
\newtheorem{notation}[theorem]{Notation}

\newcommand \bth[1] { \begin{theorem}\label{t#1} }
\newcommand \ble[1] { \begin{lemma}\label{l#1} }

\newcommand \bpr[1] { \begin{proposition}\label{p#1} }
\newcommand \bco[1] { \begin{corollary}\label{c#1} }
\newcommand \bde[1] { \begin{definition}\label{d#1}\rm }
\newcommand \bex[1] { \begin{example}\label{e#1}\rm }
\newcommand \bre[1] { \begin{remark}\label{r#1}\rm }
\newcommand \bcj[1] { \begin{conjecture}\label{j#1}\rm }

\newcommand \bnota[1] { \begin{notation}\label{n#1}\rm }
\newcommand {\eth} { \end{theorem} }
\newcommand {\ele} { \end{lemma} }

\newcommand {\epr} { \end{proposition} }
\newcommand {\eco} { \end{corollary} }
\newcommand {\ede} { \end{definition} }
\newcommand {\eex} { \end{example} }
\newcommand {\ere} { \end{remark} }
\newcommand {\ecj} { \end{conjecture} }

\newcommand {\enota} { \end{notation} }
\newcommand \thref[1]{Theorem \ref{t#1}}
\newcommand \leref[1]{Lemma \ref{l#1}}
\newcommand \prref[1]{Proposition \ref{p#1}}
\newcommand \coref[1]{Corollary \ref{c#1}}

\newcommand \deref[1]{Definition \ref{d#1}}
\newcommand \exref[1]{Example \ref{e#1}}
\newcommand \reref[1]{Remark \ref{r#1}}

\def \Rset {{\mathbb R}}         
\def \Cset {{\mathbb C}}

\def \Zset {{\mathbb Z}}

\def \Nset {{\mathbb N}}

\def \B  {{\mathcal{B}}}               
\def \AA {{\mathcal{A}}}
\def \MM {{\mathcal{M}}}

\def \QP {{\mathcal{QP}}}

\def \DD {{\mathcal{MD}}}

\def \MD {{\mathcal{MD}}}
\def \MM {{\mathcal{M}}}
\def \RD {{\mathcal{D}}}

\def \OO {{\mathcal{O}}}

\def \WW {{\mathcal{W}}}

\def \de {\delta}
\def \al {\alpha}
\def \be {\beta}

\def \la {\lambda}
\def \La {\Lambda}

\def \Om {\Omega}
\def \ga {\gamma}
\def \de {\delta}

\def \vph {\varphi}
\def \ep {\epsilon}

\def \mt  {\mapsto}
\def \ra  {\rightarrow}           
\def \lra {\longrightarrow}
\def \tra {\twoheadrightarrow}

\def \hra {\hookrightarrow}

\def \rcor {\rangle}
\def \lcor {\langle}

\def \wt {\widetilde}

\DeclareMathOperator \rank { \operatorname{rank} }

\DeclareMathOperator \norm { {\mathrm{norm}}}
\DeclareMathOperator \Ad { {\mathrm{ad}} }
\DeclareMathOperator \Rat { {\mathrm{rat}}}

\DeclareMathOperator \Span { {\mathrm{Span}} }

\DeclareMathOperator \ord { {\mathrm{ord}} }

\DeclareMathOperator \Ker { {\mathrm{Ker}} }

\DeclareMathOperator \End { {\mathrm{End}} }
\DeclareMathOperator \Gr { {\mathrm{Gr}} }
\DeclareMathOperator \Hom { {\mathrm{Hom}} }

\DeclareMathOperator \Ann { {\mathrm{Ann}} }

\DeclareMathOperator \rk {\operatorname{rank}}
\DeclareMathOperator \rad {\operatorname{rad}}

\usepackage{hyperref}
\begin{document}
\title[Adelic and Rational Grassmannians for finite dimensional algebras]
{Adelic and Rational Grassmannians for \\ finite dimensional algebras}
\dedicatory{Dedicated to Ivan Todorov on his 90th birthday with admiration and gratitude}
\author[Emil Horozov]{Emil Horozov}
\address{
Department of Mathematics and Informatics, Sofia University,
5. J. Bourchier Blvd., Sofia 1126, and, 
Institute of Mathematics and Informatics, 
Bulg. Acad. of Sci., Acad. G. Bonchev Str., Block 8, Sofia 1113, Bulgaria
}
\email{\href{horozov@fmi.uni-sofia.bg}{horozov@fmi.uni-sofia.bg}}
\author[Milen Yakimov]{Milen Yakimov}
\thanks{The research of the authors was supported by Bulgarian Science Fund grant KP-06-N62/5. The research of the second named author was also supported by NSF grant DMS-2200762.}
\address{
Department of Mathematics, Northeastern University, Boston,
MA 02115 \\
USA
}
\email{\href{m.yakimov@northeastern.edu}{m.yakimov@northeastern.edu}}
\date{}
\keywords{Wilson's adelic Grasmannian, the rational Grassmannian, finite dimensional algebras, the bispectral problem}
\subjclass[2010]{Primary 14M15; Secondary 16S32, 13E10, 37K35}
\begin{abstract} We develop a theory of Wilson's adelic Grassmannian $\Gr^{\Ad}(R)$ and Segal--Wilson's rational Grasssmannian $\Gr^{\Rat}(R)$ associated to an arbitrary finite dimensional complex algebra $R$. We provide several equivalent descriptions 
of the former in terms of the indecomposable projective modules of $R$ and its primitive idempotents, and prove that it classifies the bispectral Darboux transformations of the $R$-valued exponential function. The rational Grasssmannian $\Gr^{\Rat}(R)$ is defined by using certain free submodules of $R(z)$ and it is proved that it can be alternatively defined  via Wilson type conditions imposed in a representation 
theoretic settings. A canonical embedding  $\Gr^{\Ad}(R) \hra \Gr^{\Rat}(R)$ is constructed based on a perfect pairing between the $R$-bimodule of quasiexponentials with values in $R$ and the $R$-bimodule $R[z]$. 
\end{abstract}
\maketitle
\sectionnew{Introduction}
\label{intro}
\subsection{The classical rational and adelic Grasmannians, and bispectrality}
There are two infinite dimensional Grassmannians that play a fundamental role in the theory of the KP (Kadomtsev–Petviashvili) hierarchy \cite{vM}. The {\em{rational Grassmannian}} 
$\Gr^{\Rat}$ of Segal--Wilson 
\cite[\S 2]{SW} consists of all subspaces $M \subset \Cset(z)$ such that 
\[
h(z) \Cset[z] \subseteq M \subseteq g(z)^{-1} \Cset[z]  
\]
for some polynomials $h(z), g(z) \in \Cset[z]$ satisfying 
$\dim g(z)^{-1} \Cset[z]/M = \deg g(z)$. It parametrizes all 
solutions of the KP hierarchy that come from Krichever's 
algebro-geometric construction \cite{Kr} applied to rational curves and invertible sheaves on them. Denote by $\mathcal{C}$ the span of functionals on $\Cset[z]$ of the form $Df(z) = f^{(n)}(\la)$ for a nonnegative integer $n$ and $\la \in \Cset$. Another way of definining $\Gr^{\Rat}$, \cite[Proposition 4.6]{Wi1} is by considering the collection of subspaces 
\begin{equation}
    \label{M-def}
M= \frac{1}{g(z)} 
\{ f(z) \in \Cset[z] \mid D f(z) =0 \; \mbox{for} \; D \in V \} 
\end{equation}
where $g(z) \in \Cset[z]$ and $V$ is a finite dimensional subspace of $\mathcal{C}$ of dimension $\dim g(z)$. 

Some years after Duistermaat and Gr\"unbaum \cite{DG} posed the scalar bispectral problem, which asks for finding all analytic functions
\[
\Phi : \Om_1 \times \Om_2 \to \Cset
\]
for which there exist differential operators 
$L(x, \partial_x)$ and $\La(z, \partial_z)$ on 
$\Om_1$ and $\Om_2$ with analytic coefficients and analytic functions $\theta : \Om_1 \to \Cset$, $f : \Om_2 \to \Cset$,
such that
\begin{align}
L(x, \partial_x) \Phi(x,z) &= f(z) \Phi(x,z),  
\label{eq1}
\\
\theta(x) \Phi(x, z) &=
\La(z, \partial_z) \Phi(x,z) 
\label{eq2}
\end{align}
on $\Om_1 \times \Om_2$. Here $\Om_1, \Om_2 \subseteq \Cset$ are open, connected subsets. The scalar bispectral problem arose from questions in computer tomography and time-band limiting. There has been a great deal of research on it from the standpoints of 
integrable systems, algebraic geometry and special functions.

Wilson \cite{Wi1} made the fundamental observation that one should classify algebras of commuting bispectral differential operators instead of individual ones. He introduced the {\em{adelic Grassmannian}} 
\[
\Gr^{\Ad} \subset \Gr^{\Rat}
\]
as the subset consisting of planes $M$ given by \eqref{M-def} 
for spaces of functionals $V$ that do not mix conditions supported at different points. Wilson proved \cite{Wi1} that 
the adelic Grassmannian classifies all 
bispectral algebras of commuting differential operators of rank 1. The bispectrality of the polynomial $\tau$-functions for the KP hierarchy, which sit inside the adelic Grassmannian, was previously proved by Zubelli \cite{Z}. 
The adelic Grassmanian was related to many areas of mathematics: integrable systems of Calogero--Moser, KP, and Toda type \cite{HI,KR,Wi2}, ideal structure of the first Weyl algebra \cite{BW}, coadjoint orbits of necklace Lie algebras \cite{G}, $W$ symmetries  \cite{BN}, $A_\infty$-modules \cite{BC}, deformed preprojective algebras \cite{BCE} and others.
\medskip

The goal of this paper is to develop and 
study analogs of both the adelic and rational Grassmannians associated to every finite dimensional algebra $R$. An incomplete treatment of the adelic Grassamannian for matrix algebras was presented in \cite{Wi3}.

The paper is self-contained. In Section \ref{2.1} we gather all material on finite dimensional algebras that is needed to follow the constructions of the paper, and the reader is not expected to have a background in this subject.

\subsection{The adelic Grassmannian of a finite dimensional algebra}
\label{1.2}
One of the key methods for constructing new scalar bispectral functions from known ones is that of Darboux transformations, which goes back to \cite{DG}. In \cite{BHY2,KR} it was proved that the normalized bispectral Darboux transformations from the scalar function $\exp(xz)$ are precisely the points of Wilson's adelic Grassmannian $\Gr^{\Ad}$. Our treatment is based on the idea of performing such analysis for the  noncommutative bispectral problem with values in an arbitrary finite dimensional  complex (unital) algebra $R$. It asks for classifying the $R$-valued analytic functions
\[
\Phi : \Om_1 \times \Om_2 \to R
\]
for which there exist $R$-valued analytic differential operators 
$L(x, \partial_x)$ and $\La(z, \partial_z)$ on 
$\Om_1$ and $\Om_2$, and $R$-valued analytic functions $\theta : \Om_1 \to R$, 
$f : \Om_2 \to R$, such that
\begin{align}
L(x, \partial_x) \Phi(x,z) &= \Phi(x,z) f(z) , 
\label{eq3}
\\
\theta(x) \Phi(x, z) &=
\Phi(x,z) \La(z, \partial_z) 
\label{eq4} 
\end{align}
on $\Om_1 \times \Om_2$. The right action of $R$-valued differential operators is defined by
\begin{equation}
\label{right-act}
\Phi(x,z) \cdot  ( a(z) \partial_z^k) := (-1)^k \partial_z^k \left( \Phi(x,z) a(z) \right)
\end{equation}
and extended by linearity. In the past only the semisimple case when $R$ is the matrix algebra over $\Cset$ was considered. Some examples were constructed in \cite{BGK,BL,GI,GVZ,VZ} and other papers, while more general treatments were presented in \cite{GHY,Wi3}. In the discrete continuous setting, matrix valued bispectrality was also intensely investigated in relation to
spherical functions \cite{GPT}, matrix orthogonal polynomials \cite{C,CY}, time-band limiting \cite{CGYZ,GPZ} and other topics. 

The simplest example of an $R$-valued bispectral function is the exponential function $\exp(xz) = \exp(xz) 1$, where $1 \in R$ is the 
multiplicative identity (as it is common, we will skip writing $1$ in the future formulas for the sake of brevity).
The noncommutative bispectral Darboux transformations of it (see \thref{GenTh}) are explicitly described as the functions of the form
\[
\Phi(x,z) = P(x, \partial_x) \exp(xz),
\]
where $P(x,\partial_x)$ is a monic $R$-valued differential operator on $\Cset$ with rational coefficients which right-divides 
a monic $R$-valued differential operator $L(x, \partial_x)$ with constant coefficients. 

The set of such functions will be called the {\bf{decorated adelic Grassmannian of the algebra $R$}} and will be denoted by $\wt{\Gr}^{\Ad}(R)$. The {\bf{adelic Grassmannian}}
will be defined as the parametrizing set for the normalized bispectral $R$-valued functions. It naturally arises as a quotient of the decorated 
adelic Grassmannian $\Gr^{\Ad}(R)$.

Our first set of results is an explicit description of $\wt{\Gr}^{\Ad}(R)$ in terms of 
quiver Grassmannians associated to projective $R$-modules. We also give an explicit description of the fibers of the projection 
\[
\pi : \wt{\Gr}^{\Ad}(R) \tra \Gr^{\Ad}(R).
\]

Fix a complete set of primitive, orthogonal idempotents of the finite dimensional algebra $R$:
\[
e_1, \ldots, e_m,
\]
see \S \ref{2.1} for a review of this notion and its relations to the projective and irreducible $R$-modules. Denote by 
$\RD(R)$ the algebra of $R$-valued differential operators on $\Cset$ with rational coefficients. Consider the 
space of $R$-valued quasipolynomials 
\begin{equation}
\label{QP}
\QP(R):= \bigoplus_{\alpha \in \Cset} R[x] \exp(\alpha x).
\end{equation}
It has a canonical structure of $(\RD(R),R)$-bimodule. Recall that a finite rank, 
free (right) $R$-submodule $V$ of $\QP(R)$ is called {\em{nondegenerate}} in the sense of Etingof--Gelfand--Retakh \cite{EGR}, 
if it has a basis $F_1(x), \ldots, F_l(x) \in \QP(R)$ such that the Wronski matrix 
\[
W(F_1, \ldots, F_l) \in M_l(\QP(R))
\]
is an invertible element. We refer the reader to \S \ref{3.1} for background on quasideterminants.
\medskip
\\
\noindent
{\bf{Theorem A.}} {\em{For each finite dimensional complex {\em{(}}unital\,{\em{)}} algebra $R$, the points of the
decorated adelic Grassmannian $\wt{\Gr}^{\Ad}(R)$ are classified by the finite rank, free right $R$-submodules $V$ of the space of $R$-valued quasipolynomials 
\[
\QP(R) = \bigoplus_{\alpha \in \Cset} R[x] \exp(\alpha x)
\]
which are
\begin{enumerate}
\item[(i)] nondegenerate, and 
\item[(ii)] have bases consisting of elements of the form
\[
p_1(x) \exp(\alpha_1 x) e_1 + \cdots + p_m(x) \exp(\al_m x) e_m \quad \mbox{for some} \quad p_k(x) \in R[x], \alpha_k \in \Cset.
\]
The numbers $\al_1, \ldots, \al_l$ are not necessarily distinct and different basis elements can share the same exponents. 
\end{enumerate}
All such $R$-submodules $V$ of $\QP(R)$ are direct summands. 
The corresponding point of the decorated adelic Grassmannian is
\[
\Phi(x,z) := P(x, \partial_x) \exp(xz) \in \wt{\Gr}^{\Ad}(R)
\]
where $P(x,\partial_x) \in \RD(R)$ is the unique monic $R$-valued differential operator of order $\rank_R V$ having kernel $V$.

For a finite rank free right $R$-submodule $V$ of $\QP(R)$, the condition {\em{(ii)}} is equivalent to the condition
\begin{enumerate}
\item[(iii)] $V = \bigoplus_{k=1}^n p_k(x) \exp(\al_k x) e_{i_k} R$ for some $p_1(x), \ldots, p_n(x) \in R[x]$,
$\al_1, \ldots, \al_n \in \Cset$ and $1 \leq i_k \leq n$ {\em{(}}repetitions allowed{\em{)}}
\end{enumerate} 
and to the condition 
\begin{enumerate}
\item[(iv)] $V = \bigoplus_{\al \in \Cset} \left ( V \cap R[x] \exp(\al x) \right)$.
\end{enumerate} 
}}
\medskip

\noindent
{\bf{Remark.}} 
\begin{enumerate}
\item It is crucial that in condition (ii) in Theorem A each basis element of the module $V$ can have different exponents in front of different primitive idempotents. Requiring to have the same exponents in front of all primitive idempotents leads to very special elements of 
$\wt{\Gr}^{\Ad}(R)$. See \exref{ge2} for an illustration.
\item The nondegeneracy assumption (i) in Theorem A cannot be replaced with requiring the right $R$-submodule $V$ of $\QP(R)$ to be a direct summand. This is shown in \exref{important ex}, where $R=M_2(\Cset)$ is a semisimple algebra.
\end{enumerate}
\medskip

The meaning of the three equivalent conditions (ii)-(iv) on the right submodule $V$ of $\QP(R)$ is that (ii) provides a 
basis of $V$, (iii) describes the decomposition of $V$ into indecomposable projectives, and (iv) describes $V$ with respect to the 
direct sum decomposition \eqref{QP} of $\QP(R)$. All three are adelic type conditions in the sense of Wilson \cite{Wi1,Wi2}. 
In condition (iii), $n =  m \rank_R V = m \ord P$.

A (right) submodule $V$ of $\QP(R)$ is a direct summand of $\QP(R)$  if and only if it is a direct summand of one (and thus of any) submodule 
of $\QP(R)$ of the form 
\[
X= (R + \cdots + R x^{N_1-1}) \exp{\al_1 x} \bigoplus \cdots \bigoplus (R + \cdots + R x^{N_k-1}) \exp(\al_k x)
\]
such that $V$ lies in $X$. This follows from the fact that, if $V \oplus W =U$ for modules $V,W,U$ and $X$ is a submodule of $U$ containing $V$, 
then $V$ is a direct summand of $X$, namely $X = V \oplus (X \cap W)$.

In the case when $R= \Cset$, the condition (i) is automatically satisfied but is nontrivial apart from this case.

Define the normalization of an $R$-valued bispectral function 
\[
\Phi(x,z) = P(x, \partial_x) \exp(xz) \in  \wt{\Gr}^{\Ad}(R)
\]
corresponding (via Theorem A) to the right $R$-submodule $V$ of $\QP(R)$ with basis 
\[
p_{1j}(x) \exp(\alpha_{1j} x) e_1 + \cdots + p_{mj}(x) \exp(\al_{mj} x) e_m \quad \mbox{for} \quad p_{kj}(x) \in R[x], \alpha_{kj} \in \Cset, 1 \leq j \leq l
\]
to be the function
\begin{equation}
\label{normal}
\Phi_{\norm}(x,z) := \Big[  P(x, \partial_x) \exp(xz) \Big] . \Big[ \prod_j(z- \al_{1j}) e_1 + \cdots +  \prod_j(z- \al_{mj}) e_m \Big]^{-1}.
\end{equation}
The set of such normalized $R$-valued bispectral functions will be called the {\bf{adelic Grassmannian of the algebra $R$}} and will be denoted by $\Gr^{\Ad}(R)$. 

Next we turn our attention to the fibers of the canonical projection 
\[
\pi : \wt{\Gr}^{\Ad}(R) \tra \Gr^{\Ad}(R), \quad \Phi(x,z) \mapsto \Phi_{\norm}(x,z).
\]
Consider the partial order on $\wt{\Gr}^{\Ad}(R)$ given by 
\[
\Phi(x,z) = P(x, \partial_x) \exp(xz) \prec \Psi(x,z) = S(x, \partial_x) \exp(xz) \in \wt{\Gr}^{\Ad}(R)
\]
whenever
\[
S(x, \partial_x) = P(x, \partial_x)  \prod_j \big( \partial_x- \gamma_{1j} e_1 - \cdots -  \gamma_{mj} e_m \big)
\]
for some scalars $\gamma_{kj} \in \Cset$, $1 \leq k \leq m$, $1 \leq j \leq l$.
\medskip

\noindent
{\bf{Theorem B.}}
\begin{enumerate}
\item {\em{The fiber of $\pi : \wt{\Gr}^{\Ad}(R) \tra \Gr^{\Ad}(R)$ containing  $\Phi(x,z) \in \wt{\Gr}^{\Ad}(R)$ consists of all 
$\Psi(x,z) \in \wt{\Gr}^{\Ad}(R)$ such that $\Phi(x,z)$ and $\Psi(x,z)$ have a common successor.}}
\item {\em{Let, under the classification of Theorem A, $\Phi(x,z) \in \wt{\Gr}^{\Ad}(R)$
corresponds to the right $R$-submodule
\[
V := \bigoplus_{k=1}^n p_k(x) \exp(\al_k x) e_{i_k} R \subset \QP(R)
\] 
of $\QP(R)$.
The immediate successors of $\Phi(x,z)$ in $\wt{\Gr}^{\Ad}(R)$ are parametrized by $\Cset^m$ 
and correspond to the right $R$-submodules
\[
V_{(\ga_1, \ldots \ga_n)}:= \big( \bigoplus_{i=1}^m \exp(\ga_i x) e_i R \Big) \bigoplus \Big( \bigoplus_{k=1}^n F_{A_{\al_k}} (p_k(x) e_{i_k}) \exp(\al_k x)  R \Big)
\subset \QP(R),
\]
for $(\ga_1, \ldots, \ga_m) \in \Cset^m$,
where $A_{\al} \in \End_{\Cset} (R)$ is the operator of left multiplication by $(\al - \ga_1) e_1 + \cdots + (\al- \ga_m) e_m$.}}
\end{enumerate}
\subsection{The rational Grassmannian of a finite dimensional algebra}
\label{1.3}
Now consider $R[z]$ as a left(!) $R$-module.  
Define {\bf{the rational Grassmannian}} $\Gr^{\Rat}(R)$ to be the set of free (left) $R$-submodules $M \subset R(z)$ satisfying the following two conditions 
\begin{enumerate}
\item[(i)] $h(z) R[z] \subseteq M \subseteq g(z)^{-1} R[z]$ for some $h(z), g(z) \in \Cset[z]\backslash \{0\}$ and
\item[(ii)]  $M/ \left( h(z) R[z] \right)$ is a direct summand of $g(z)^{-1} R[z] / \left( h(z) R[z]\right)$ as a left $R$-module
such that
\[
\rk_R \big( g^{-1}(z) R[z] / M \big) = \deg g(z).
\]
\end{enumerate}

In the scalar case $R = \Cset$ the direct summand assumption in condition (ii) is automatically satisfied and we recover the classical Segal--Wilson adelic Grassmannian due to \cite[Proposition 4.6]{Wi1}. 


Define the pairing
\begin{equation}
\label{form}
\lcor .,. \rcor :  R[z] \times \QP(R) \to R
\end{equation}
by
\begin{equation}
\label{pairing-formula}
\lcor p(z),  x^n \exp(\al x) r   \rcor := p^{(n)}(\al) r \quad \forall \al \in \Cset, n \in \Nset, \; \; p(z) \in R[z].
\end{equation}
It is left $R$-linear in the first argument and right $R$-linear in the second
\begin{equation}
\label{linearity}
 \lcor r p, f \rcor = r \lcor p, f \rcor, \quad \lcor p, f r \rcor = \lcor p, f \rcor r \quad \forall r \in R, f \in \QP(R), p \in R[z], 
\end{equation}
and satisfies 
\begin{equation}
\label{mixed-invar}
\lcor p r, f \rcor = \lcor p, r f\rcor \quad \forall r \in R, f \in \QP(R), p \in R[z].
\end{equation}
For a right $R$-submodule $V \subset \QP(R)$, consider its orthogonal complement 
\[
V^\perp \subset R[z],
\]
which by \eqref{linearity} is a left $R$-submodule. 

Next we define a map 
\[
\wt{\iota} : \wt{\Gr}^{\Ad}(R) \to \Gr^{\Rat}(R).
\]
For an element $\Phi(x,z) \in \Gr^{\Ad}(R)$, which under the classification of Theorem A corresponds to $V \subset \QP(R)$ 
with basis
\[
p_{1j}(x) \exp(\alpha_{1j} x) e_1 + \cdots + p_{mj}(x) \exp(\al_{mj} x) e_m \quad \mbox{for} \quad p_{kj}(x) \in R[x], \alpha_{kj} \in \Cset, 1 \leq j \leq l,
\]
set
\begin{align*}
\wt{\iota}(\Phi) &:= V^\perp . \Big[ \prod_j(z- \al_{1j}) e_1 + \cdots +  \prod_j(z- \al_{mj}) e_m \Big]^{-1} 
\nn
\\
&= V^\perp . \Big[ \frac{\theta(z)}{\prod_j(z- \al_{1j})} e_1 + \cdots +  \frac{\theta(z)}{\prod_j(z- \al_{mj})} e_m \Big] . \theta(z)^{-1},
\end{align*}
where
\[
\theta(z) = \prod_{ij} (z - \al_{ij}).
\]
In Section \ref{5.3} we show that $\wt{\iota}$ takes values in $\Gr^{\Rat}(R)$. Our third main result is that $\wt{\iota}$ descends to an embedding $\Gr^{\Ad}(R) \hra \Gr^{\Rat}(R)$.
\medskip
\\
\noindent
{\bf{Theorem C.}} {\em{The map $\wt{\iota} : \wt{\Gr}^{\Ad}(R) \to \Gr^{\Rat}(R)$ descends under the projection \\ $\pi : \wt{\Gr}^{\Ad}(R) \tra \Gr^{\Ad}(R)$ to an embedding
$\iota : \Gr^{\Ad}(R) \hra \Gr^{\Rat}(R)$:
\[
\begin{tikzcd}[column sep=small]
&  \wt{\Gr}^{\Ad}(R) \arrow[two heads, "\pi" above]{dl}
\arrow{dr}{\wt{\iota}} & \\
\Gr^{\Ad}(R) \arrow[hookrightarrow]{rr}{\iota}& & \Gr^{\Rat}(R)
\end{tikzcd}
\]
}}
\medskip

Combining Theorems B and C, for every finite dimensional algebra $R$, we obtain the chain of maps 
\begin{equation}
    \label{3Grassm}
\wt{\Gr}^{\Ad}(R) \tra \Gr^{\Ad}(R) \hra \Gr^{\Rat}(R).
\end{equation}
\subsection{Organization of the paper}
The paper is organized as follows.
Section \ref{GenTh} contains background material on 
finite dimensional algebras and bispectral Darboux transformations, as well as results on the decorated adelic Grassmannian. The first part of Section \ref{Facto} contains background material on qiasideterminants and factorizations of 
differential operators with values in a differential algebra following \cite{EGR,GR1,GR2,K}, and the second part treats factorization of $R$-valued differential operators with meromorphic coefficients. 
Section \ref{class} contains the proofs of Theorems A and B. Section \ref{rational} develops the rational Grassmannian of a finite dimensional algebra $R$ and proves Theorem C. 
Section \ref{ex} contains examples illustrating the constructions of the Grassmannains $\wt{\Gr}^{\Ad}(R)$, $\Gr^{\Ad}(R)$, and $\Gr^{\Rat}(R)$ and the maps between them \eqref{3Grassm}. 

It is natural to ask how to equip the three Grasmannians $\wt{\Gr}^{\Ad}(R)$, $\Gr^{\Ad}(R)$, and $\Gr^{\Rat}(R)$ with ind-algebraic structures in a compatible fashion with the maps \eqref{3Grassm} and what the functoriality properties of the constructions are. We will address this in a later publication. We will also investigate 
partitions of $\Gr^{\Ad}(R)$ into generalized 
Calogero--Moser spaces and embeddings of the latter as coadjoint orbits of necklace Lie algebras \cite{G} associated to $R$. 
\subsection{Terminology and notation}
All rings and algebras considered in this paper will be assumed to be unital, and by 
an algebra we mean an associative algebra. We denote
\[
\Zset_+ := \{1, 2, \ldots\} \quad \mbox{and} \quad \Nset:= \{0,1,\ldots\}.
\]
\medskip

\noindent
{\bf Acknowledgements.} M.Y. thanks the hospitality of the Institute of Mathematics and Informatics of the Bulgarian Academy of Sciences where parts of this project were completed. 
\sectionnew{Finite dimensional algebras and noncommutative bispectral Darboux transformations}
\label{GenTh}
In this section we review background material on finite dimensional algebras, noncommutative localization and the construction
of bispectral functions via noncommutative Darboux transformations. We then describe the regular elements of 
the algebra $R[x]$ for a finite dimensional algebra $R$
and use this fact to describe the bispectral Darboux transformations of the 
exponential function in an $R$-valued setting and to define the decorated adelic Grassmannian of $R$. 
\subsection{Finite dimensional algebras}
\label{2.1}
To facilitate the reading of the paper, here we gather material about finite dimensional algebras; 
we will restrict to the case of finite dimensional modules, which is what is needed for the paper. 
The section contains all details needed to follow the paper. If the reader prefers to have additional information, we refer to Chapters I and III of \cite{ASS}.  

The identity element $1 \in R$ has a {\em{decomposition}} 
\begin{equation}
\label{idemp-decomp}
1 = e_1 + \ldots + e_m,
\end{equation}
where $\{e_1, \ldots, e_m\}$ is a {\em{complete set of primitive, orthogonal idempotents}}, meaning: 
\begin{enumerate}
\item[(i)] $e_i \in R$ are {\em{idempotents}} ($e_i^2 = e_i$ for all $1 \leq i \leq m$), which are
\item[(ii)] {\em{pairwise orthogonal}} ($e_i e_j = 0$ for all $i \neq j$) and
\item[(iii)] {\em{indecomposable}} ($e_i \neq e + e'$ for two nonzero orthogonal idempotents $e$ and $e'$),
\end{enumerate}   
\cite[Section I.4]{ASS}. 

An $R$-module $P$ is {\em{projective}} if it possesses the {\em{lifting property}} that for all surjective homomorphisms 
\[
\phi : M \to N,
\]
a homomorphism $f : P \to N$ can be lifted to a 
homomorphism $f : P \to M$ such that $f = \phi g$.
An $R$-module is projective if and only if it is a {\em{direct summand of a free module}}.  
An $M$-module is {\em{indecomposable}} if $M \not{\cong} M' \oplus M''$ for two proper submodules $M', M'' \subset M$, 
and simple if it has no nontrivial proper submodules. 

By the {\em{Krull--Schmidt property}}, \cite[Theorem 4.10]{ASS}, every two decompositions of a module into direct summands have isomorphic collections 
of direct summands (counted with multiplicities).

In terms of the decomposition \eqref{idemp-decomp}, the isomorphism classes of {\em{indecomposable projective right $R$-modules}} are the modules
\[
e_1 R, \ldots, e_m R,
\]
\cite[Theorem 5.8]{ASS}. However, note that some of them can be isomorphic, see \exref{ex-fin-dim-alg}(i) below.  
 
The {\em{(Jacobson) radical}} $\rad R$ of the algebra $R$ is the intersection of all maximal left ideals (by definition, the latter are proper), or equivalently the intersection of all maximal right ideals, \cite[Definition 1.2]{ASS}. It is 
characterized by
\[
\rad R = \{ r \in R \mid 1 - r t \; \; \mbox{is an invertible element of $R$ for all $t \in R$} \}, 
\]
\cite[Lemma 1.3]{ASS}
and is {\em{nilpotent}}
\[
(\rad R)^\ell =0 \quad \mbox{for some positive integer $\ell$},
\]
\cite[Corollary 2.3]{ASS}. 
\bex{ex-fin-dim-alg} 
(i) The matrix algebra $M_n(\Cset)$ is semisimple and  
\[
1 = e_{11} + \ldots + e_{nn}
\]
is an indecomposable orthogonal decomposition of the identity, where $e_{ij}$ denote the standard matrix entries. 
The indecomposable projective modules $e_{ii} R$ are all isomorphic to the standard module $\Cset^n$ (realized as row-vectors under the right action of $M_n(\Cset)$), which is the only irreducible $R$-module. 

An algebra $R$ is called {\em{semisimple}} if all of its finite dimensional modules are direct sums of simple ones. 
By the {\em{Wedderburn--Artin theorem}} \cite[Theorem 3.4]{ASS} this is equivalent to the condition that its Jacobson radical is 0, and the finite dimensional semisimple algebras over $\Cset$ are directs products of matrix algebras over $\Cset$.  

(ii) The {\em{algebra of dual numbers}} is the two-dimensional algebra 
\[
R:=\Cset[\epsilon]/(\epsilon^2).
\]
This algebra is not semisimple. In fact, it is a {\em{local algebra}}:  
\[
\rad R = \Cset \ep
\]
and all elements outside $\rad R$ are invertible. Moreover, 
$1$ is a primitive idempotent and $R$ is an idecomposable projective module (which is unique 
up to isomorphism). 
\eex
For additional interesting examples of nonsemisimple finite dimensional algebras, we refer the reader to Section \ref{ex}. 
\subsection{The generalized bispectral setting and noncommutative bispectral Darboux transformations}
\label{2.2}
A {\em{multiplicative subset}} $E$ of a ring $A$ is a subset which is closed under multiplication and contains the identity element of $A$. 
A multiplicative subset $E$ of $A$ 
consisting of {\em{regular elements}} (i.e., non-zero divisors) is called a (left) {\em{Ore set}} if for all $s \in E$ and $a \in A$,
there exist $t \in E$ and $b \in A$ such that $t a = b s$. In this case one can form the {\em{quotient ring}} $A[E^{-1}]$
consisting of functions of the form $s^{-1} a$ for $s \in E$, $a \in A$ with a natural equivalence relation, and addition and multiplication 
operations, see \cite[pp. 105-110]{GW} for details. There is a canonical embedding $A \hra A[E^{-1}]$ given by $a \mt 1^{-1} a$. 
Denote by
\[
r(A) \; \; \mbox{the set of all regular elements of $A$},
\]
i.e., elements that are not left or right 0 divisors. 
It is a multiplicative subset of $A$.

Next, we describe the setting of noncommutative Darboux transformations and Theorem 2.1 from \cite{GHY}
showing that they produce bispectral functions; the commutative case appeared in \cite[Theorem 4.2]{BHY1}. 
Consider the following: 
\medskip

\noindent
{\bf{Generalized bispectral setting}}:

(1) Let $(B_1, B_2)$ be a pair of algebras and $b : B_1 \stackrel{\cong}{\lra} B_2$
be an isomorphism. 

(2) Let $M$ be a $(B_1,B_2)$-bimodule and $\psi \in M$ be such that 
\[
Q \psi = \psi b(Q), \quad \forall Q \in B_1.
\] 
We call $B_1$ and $B_2$ the {\em{left}} and {\em{right Fourier algebras}} of $\psi$ and $b :  B_1 \stackrel{\cong}{\lra} B_2$
the {\em{Fourier isomorphism}}. 

(3) Let $A_i$ and $K_i$ be subalgebras of $B_i$ such that 
\[
b(A_1) = K_2, \quad b(K_1) = A_2
\]
and which satisfy the following conditions:
\begin{enumerate}
\item[(i)]
$r(K_i)$ are Ore subsets of $B_i$ consisting 
of regular elements and

\item[(ii)] $\Ann_f(M) = 0$ for all $f \in r(K_i)$.  
\end{enumerate}
The algebras $A_1$ and $A_2$ will be called the {\em{left}} and {\em{right spectral algebras}} of $\psi \in M$.
\medskip

We think of $K_i$ as generalizations of the algebras of functions in $x$ and $z$ in \eqref{eq3}--\eqref{eq4} 
and of $A_i$ as generlizations of the algebras of differential operators in $x$ and $z$. In this way, the properties
\begin{align*}
L \psi &= \psi b(L) , 
\\
b^{-1}(\La) \psi &=
\psi \La
\end{align*}
for nonzero elements $L \in A_1$, $\La \in A_2$, represent a generalized bispectrility property of $\psi$. 

The algebras
\[
\wt{B}_i := B_i[r(K_i)^{-1}]
\]
will be called the {\em{left}} and {\em{right localized Fourier algebras}} of $\psi \in M$. Conditions (i) and (ii) ensure that the left/right 
Fourier algebras $B_i$ embed in their localized counterparts and that the $(B_1, B_2)$-bimodule
$M$ embeds in the $(\wt{B}_1,\wt{B}_2)$-bimodule $M[r(K_1)^{-1}][r(K_2)^{-1}]$,
see e.g. Theorem 10.8 in \cite{GW} and its proof. 

\bex{Fourier} Consider the  first Weyl algebras in the variables $x$ and $z$
\[
B_1 := \Cset [x, \partial_x] \quad \mbox{and} \quad B_2 := \Cset [z, \partial_z],
\]
and the isomorphism
\[
b : B_1 \stackrel{\cong}{\lra} B_2, 
\quad \mbox{given by} \quad 
b(\partial_x) = z, \; \; b(x) = - \partial_z.
\]
Up to a sign, the latter is the classical Fourier map on differential operators. 
Let $M$ be the space of meromorphic functions on $\Cset \times \Cset$
and $\psi:= \exp(x,z) \in M$. The space $M$ is a $(B_1,B_2)$-bimodule structure
by using the standard left action of the Weyl algebra $B_1$ on $M$ 
and the right action \eqref{right-act} of the Weyl algebra $B_2$ on $M$.

Choose 
\[
K_1 := \Cset[x] \subset \Cset [x, \partial_x] \quad \mbox{and} \quad
K_2 := \Cset[z] \subset \Cset [z, \partial_z]
\]
be the subalgebras of the Fourier algebras $B_i$ consisting of differential operators 
of orders 0. The corresponding left and right spectral algebras of $\psi = \exp(xz)$ are
\[
A_1 = b^{-1}(K_2) = \Cset [\partial_x] \quad \mbox{and} \quad
A_2 = b(K_1) = \Cset [\partial_z]. 
\]
They describe all differential operators for which $\psi = \exp(xz)$ is a bispectral function, i.e., 
\begin{align*}
L(x, \partial_x) \psi &= \psi f(z) , 
\\
\theta(x) \psi  &=
\psi \La(z, \partial_z) 
\end{align*}
if and only if $L(x, \partial_x) \in A_1$, $\La(z, \partial_z) \in A_2$ and $f = b(L) \in K_1$, $\theta = b^{-1}(\La) \in K_1$. 

The multiplicative subsets of regular elements of $K_1$ and $K_2$ are 
\begin{equation}
\label{rK}
r(K_1) = \Cset[x] \backslash \{ 0 \} \quad \mbox{and} \quad r(K_2) = \Cset[z] \backslash \{ 0 \},
\end{equation}
and the localized left and right Fourier algebras are the algebras
\[
\wt{B}_1 = \Cset(x) [\partial_x] \quad \mbox{and} \quad \wt{B}_2 = \Cset(z) [\partial_z]
\]  
of differential operators with rational coefficients in the variables $x$ and $z$. 
\qed
\eex
In the above generalized bispectral setting, one would like to 
construct new bispectral functions in the module $M[r(K_1)^{-1}][r(K_2)^{-1}]$
out of the initial bispectral function $\psi \in M$. This is achieved by the next theorem.

\bth{GenTh} \cite[Theorem 2.7]{CGYZ}
Assume the generalized bispectral setting. 

If $L = P'  Q'$ with $L \in A_1$ and $P', Q' \in \wt{B}_1$ with $P'$ that is not a left zero divisor of $\wt{B}_1$, then 
\[
Q' \psi \in M[r(K_1)^{-1}][r(K_2)^{-1}]
\] 
is a bispectral function. 

More precisely, if $L = S g^{-1} Q$ for some $L \in A_1$, $S,Q \in B_1$, and $g \in r(K_1)$, then 
$\vph := Q \psi = \psi b(Q) \in M[r(K_1)^{-1}][r(K_2)^{-1}]$ satisfies
\begin{align}
\label{thm-eq1}
\left( QS g^{-1} \right) \vph &= \vph b(L), \\
g \vph &= \vph \left( b(L)^{-1} b(S) b(Q) \right)
\label{thm-eq2}
\end{align} 
{\em{(}}note that $b(L) \in K_2${\em{)}}. 
\eth
\subsection{Regular elements of $R[x]$}
\label{2.3}
For the rest of the paper, let $R$ be a finite dimensional algebra over $\Cset$, as in the beginning of the section. 
In order to deal with the localizations in \thref{GenTh} in the case of $R$-valued differential operators with polynomial coefficients, 
we will need a characterization of the regular elements $R[x]$. 
The characterization of that set is obtained next. 

For $p(x) \in R[x]$, denote the $\Cset$-linear operator
\[
[ p(x) \cdot] : R \to R[x],  \quad \mbox{given by} \quad [p(x) \cdot](a) :=  p(x) a.  
\]
Evaluating at different $x \in \Cset$, we obtain $\Cset$-linear endomorphisms of $R$ whose determinants will be denoted by $\det [p(x) \cdot] \in \Cset[x]$. 
Analogously, $\det[ \cdot p(x)]$ denotes the determinants of the $\Cset$-linear endomorphisms of $R$ 
given by right multiplication by $p(x)$ for $x \in \Cset$.
\ble{det} The following are equivalent for $p(x) \in R[x]$ for a finite dimensional algebra
\begin{enumerate}
\item $p(x)$ is not a left zero divisor,
\item $p(x)$ is not a right zero divisor,
\item $\det [p(x) \cdot] \in \Cset[x]$ is a non-zero polynomial,
\item $\det [\cdot p(x)] \in \Cset[x]$ is a non-zero polynomial,
\item there exists $q(x) \in R[x]$ such that $p(x) q(x) = q(x) p(x) = s(x) \in \Cset[x] \backslash \{ 0 \}$.  
\end{enumerate}
\ele
\begin{proof} We first prove that (1) $\Leftrightarrow$ (3) $\Leftrightarrow$ (5') $p(x) q(x) = s(x) \in \Cset[x] \backslash \{ 0 \}$
for some $q(x) \in R[x]$. 

Clearly (5') $\Rightarrow$ (3). For (3) $\Rightarrow$ (1), note that, if $p(x) f(x) = 0$ for some non-zero 
$f(x) \in R[x]$, then $\det[ p(x_0) \cdot ] =0$ for all $x_0 \in \Cset$ that are not roots of $f(x)$, so $\det[ p(x) \cdot ] =0$. 

For (1) $\Rightarrow$ (5'), assume that $\det[ p(x) \cdot ] \neq 0$. For all $x_0 \in \Cset$ such that $\det [ p(x_0) \cdot ] \neq 0$, 
$[p(x_0) \cdot] : R \to R$ is an invertible $\Cset$-linear operator. Applying the Kramer's rule to it, we obtain that 
there exists $q(x) \in R[x]$ such that 
\[
p(x) q(x) = \det [p(x) \cdot] \in \Cset[x] \backslash \{ 0 \}.
\]

Analogously, one shows that (2) $\Leftrightarrow$ (4) $\Leftrightarrow$ (5'') $q'(x) p(x) = s'(x) \in \Cset[x] \backslash \{ 0 \}$
for some $q'(x) \in R[x]$. 

For (5') $\Rightarrow$ (2) note that, if $p(x) q(x) = s(x) \in \Cset[x] \backslash \{ 0 \}$
for some $q(x) \in R[x]$, then $f(x) p(x) =0$ for $f(x) \in R[x]$, implies $f(x) s(x) = f(x) p(x) q(x) =0$, and so, 
$f(x)=0$ because $\Cset[x] \backslash \{ 0 \}$ consists of regular elements of $R[x]$. 

Analogously, one shows that (5'') $\Rightarrow$ (1). Therefore (1), (2), (3), (4), (5') and (5'') are equivalent.

It remains to show that (5) $\Leftrightarrow$ (5')$\&$(5''). Obviously, (5) $\Rightarrow$ (5')$\&$(5''). For the opposite implication, 
assume that (5') and (5'') are satisfied. Thus, $p(x) q(x), q'(x) p(x) \in \Cset[x] \backslash \{0\}$ for some 
$q(x), q'(x) \in R[x]$. After multiplying the two products with terms in $\Cset[x]$ we can assume that 
\[
p(x) q(x) = q'(x) p(x) = s(x) \in \Cset[x] \backslash \{0\}.
\] 
Therefore, 
\[
q'(x) s(x) = q'(x) p(x) q(x) = s(x) q(x),  
\]
which implies that $q(x) = q'(x)$ because the elements of $\Cset[x] \backslash \{ 0 \}$ are regular and belong to the center of $R[x]$. 
\end{proof}

Considering $r \in R \subset R[x]$ in \leref{det}, one obtains the following 
corollary. In it we denote by $(r \cdot)$ and $(\cdot r)$ the $\Cset$-linear endomorphisms of $R$ 
given by left and right multiplication by $r$, and by 
\begin{equation}
\label{determinants}
\det [r \cdot],  \det [\cdot r] \in \Cset
\end{equation}
their determinants. 
\bco{regR} The following are equivalent for an element $r$ of a finite dimensional algebra $R$:
\begin{enumerate}
\item $r$ is not a left zero divisor,
\item $r$ is not a right zero divisor,
\item $\det [r \cdot] \neq 0$,
\item $\det [\cdot r] \neq 0$,
\item $r$ is invertible. 
\end{enumerate}
\eco
\subsection{Definition of the decorated adelic Grassmannian}
\label{2.4}
The adelic Grassmannian of $R$ studied in this paper arises from a classification of the 
$R$-valued bispectral functions as in \eqref{eq3}--\eqref{eq4} that come from 
\thref{GenTh}. This has to do with the following setting that generalizes \exref{Fourier}. Let
\[
B_1 := R [x, \partial_x] \quad \mbox{and} \quad B_2 := R[z, \partial_z]
\]
be the first Weyl algebras in the variables $x$ and $z$ with values in $R$, and $b : B_1 \stackrel{\cong}{\lra} B_2$ be the isomorphism, given by
\[
b(\partial_x) = z, b(x) = - \partial_z, b(r) = r, \; \forall r \in R.
\]
Let 
\begin{equation}
\label{mer-fn}
M := \MM(R)
\end{equation}
be the space of $R$-valued meromorphic functions on $\Cset \times \Cset$. It has a canonical 
$(B_1,B_2)$-bimodule structure, where $B_1$ acts on $\MM(R)$ via the standard left action and $B_2$ via the right action \eqref{right-act}.
Denote 
\[
\psi:= \exp(x,z) \in \MM(R).
\]
We have that 
\begin{equation}
\label{phi}
S \psi = \psi b (S), \quad \forall S \in B_1 = R [x, \partial_x].
\end{equation}
Denote the subalgebras $K_i \subset B_i$ consisting of functions (i.e., operators of degree 0),
\[
K_1 := R[x], \quad K_2 := R[z].
\]
Their images under $b$ and $b^{-1}$ are the subalgebras of operators with constant coefficients: 
\[
A_1 := R [\partial_x], \quad A_2 := R [\partial_z]. 
\]
As in \exref{Fourier}, 
they describe all differential operators for which $\psi = \exp(xz)$ is a bispectral function:
\begin{align*}
L(x, \partial_x) \psi &= \psi f(z) , 
\\
\theta(x) \psi  &=
\psi \La(z, \partial_z) 
\end{align*}
if and only if $L(x, \partial_x) \in A_1$, $\La(z, \partial_z) \in A_2$ and $f = b(L) \in K_1$, $\theta = b^{-1}(\La) \in K_1$. 

The multiplicative subsets $r(K_i) \subset B_i$ 
are characterized in \leref{det}. The equivalent condition (5) in the lemma implies that they are Ore subsets and
\begin{equation}
\label{locB}
B_1 [ r(K_1)^{-1}] \cong B_2 [ r(K_2)^{-1}] \cong \RD(R). 
\end{equation}
(Recall from \S \ref{1.2} that $\RD(R)$ denotes the algebra of $R$-valued differential operators on $\Cset$ with rational coefficients.)

Since we have an embedding $B_1 = R[x, \partial_x] \subset \RD(R)$, in order to prove \eqref{locB}, all we need to show is that the elements 
$r(R[x])$ are invertible in $\RD(R)$ and that $\RD(R)$ consists of the left fractions $b^{-1} a$ (as well as the right fractions 
$a b^{-1}$) with $a \in R[x, \partial_x]$ and $b \in r(R[x])$, see \cite[Theorem 6.2]{GW}. This follows from \leref{det}.
This verifies condition (i) in the definition of generalized bispectral setting in \S \ref{2.2} and describes the 
left and right localized Fourier algebras of $\psi = \exp(xz) \in M$. The equivalent condition (5) in \leref{det} also implies that 
$r(R[x])$ and $r(R[z])$ are units (invertible elements) of $\MM(R)$, which verifies condition (ii) 
in the definition of generalized bispectral setting in Section \ref{2.2} and shows that 
\[
\MM(R)[r(R[x])^{-1}][r(R[z])^{-1}] = \MM(R).
\] 
Finally, we note that every monic differential operator in $\RD(R)$ is not a left or right zero divisor of $\RD(R)$. 
We obtain the following corollary of \thref{GenTh}:
\bco{bisp} Let $R$ be a finite dimensional complex algebra. All functions of the form
\[
\Phi(x,z) = P(x, \partial_x) \exp(xz) \in \MM(R),
\]
where $P(x,\partial_x)$ is a monic $R$-valued differential operator on $\Cset$ with rational coefficients which right-divides 
a monic $R$-valued differential operator $L(x, \partial_x)$ with constant coefficients, are bispectral. 

More precisely, if \[
L(x, \partial_x) = Q'(x, \partial_x) h(x)^{-1} g(x)^{-1} P'(x, \partial_x)
\]
for $L(x, \partial_x) \in R[\partial_x]$, $P'(x, \partial_x), Q'(x, \partial_x) \in R[x,\partial_x]$, $g(x), h(x) \in r(R[x])$, the meromorphic 
function
\[
\Phi(x,z) := g(x)^{-1} P'(x,\partial_x) \Psi(x,z) : \Cset \times \Cset \to R
\]
satisfies
\begin{align*}
\left( g(x)^{-1} P'(x, \partial_x) h(x)^{-1} Q'(x, \partial_x) \right) \Phi(x,z) 
&= \Phi(x,z) b(L)(z) \\
g(x) h(x) \Phi(x,z) &= \Phi(x,z) \left( b(L)(z)^{-1} b(Q')(z, \partial_z) b(P')(z, \partial_z) \right).
\end{align*} 
\eco 
\bde{adGr} We define the {\bf{decorated adelic Grassmannian}} of the finite dimensional algebra $R$ to be the set $\wt{\Gr}^{\Ad}(R)$ consisting of 
$R$-valued meromorphic functions on $\Cset \times \Cset$ of the form:
\[
\Phi(x,z) = P(x, \partial_x) \exp(xz) \in \MM(R),
\]
where $P(x,\partial_x) \in \RD(R)$ is a monic $R$-valued differential operator on $\Cset$ with rational coefficients which right-divides 
a monic $R$-valued differential operator $L(x, \partial_x)$ with constant coefficients. Equivalently, $\wt{\Gr}^{\Ad}(R)$ parametrizes 
the set of differential operators $P(x, \partial_x) \in \RD(R)$ with this property.
\ede
\bre{big} Although we start with the bispectral function $\exp(xz)$ that takes values in $\Cset \subset R$,
the fact that all Darboux transformations are performed in an $R$-valued setting makes 
$\wt{\Gr}^{\Ad}(R)$ much bigger and more complicated than $\wt{\Gr}^{\Ad}(\Cset)$.
\ere
\sectionnew{Factorization of $R$-valued meromorphic differential operators}
\label{Facto}
In the first part of this section we gather background material on qiasideterminants and factorizations of 
differential operators with values in a differential algebra. We then investigate in detail the factorizations of $R$-valued differential operators with meromorphic coefficients in 
terms of what right $R$-submodules of their kernels 
can occur in such factorizations. Here, the assumption that the target algebra $R$ is finite dimensional is essential. 
\subsection{Results in the differential algebra setting}
\label{3.1}
First we prove a factorization theorem for certain kinds of differential operators in the setting of differential rings considered in \cite{EGR}.

We briefly recall the notion of {\em{quasideterminants}} of Gelfand and Retakh \cite{GR1,GR2}.
Let $A$ be a unital ring and $Y= (y_{ij}) \in M_n(A)$. For $i,j \in [1,n]$, the $i$-th row and $j$-th column 
of $Y$ will be denoted by $r_i(Y)$ and $c_j(Y)$. The matrix obtained by removing the $i$-th row and the $j$-th column of $Y$
will be denoted by $Y^{ij}$. For row or column vectors $r$ and $c$, $r^{(i)}$ and $c^{(i)}$ will denote the vectors 
obtained by removing their $i$-th entries. If $Y^{ij} \in M_{n \times n}(A)$ is invertible, the quasideterminant  $|Y|_{ij}$ 
is defined by 
\[
|Y|_{ij}   :=    y_{ij} - r_i(Y)^{(j)} (Y^{ij})^{-1} c_j(Y)^{(i)} \in A
\]
when the inverse $(Y^{ij})^{-1}$ exists. 
This is a recursive definition because the inverse of a square matrix with entries in $A$ is simply the 
matrix whose entries are the inverses of the quasideterminants of the original matrix, 
see \cite[Eq. (1.6)]{GR1} and \cite[Theorem 1.2.1]{GR2}. 

Let $A$ be a unital differential ring with differential $\delta$ and $A[\partial]$
be the corresponding ring of differential operators with coefficients in $A$.
In this setting, $A$ is a left $A[\partial]$-module by $\partial^n F := \delta^n F$.
The {\em{Wronski matrix}} of $F_1, \ldots, F_n \in A$ is defined by 
\begin{equation}
\label{Wr}
W(F_1, \ldots, F_n) :=   \left(
\begin{matrix}
F_1 & \ldots & F_n \\ 
\delta F_1 & \ldots & \delta F_n  \\ 
\ldots & \ldots & \ldots \\ 
\delta^{n-1} F_1 & \ldots &\delta^{n-1} F_n
\end{matrix} \right). 
\end{equation}
\bde{nondeg0} \cite{EGR}
A collection of elements $\{F_1, \ldots, F_n\}$ of $A$ is called {\em{nondegerate}} if $W(F_1, \ldots, F_n)$ is 
an invertible element of the algebra of $n \times n$-matrices with entries in $A$. 
\ede
The next theorem describes the monic differential operators in $A[\partial]$ 
annihilating nondegenerate collections of elements of $A$. 
\bth{EGR} \cite[Theorem 1.1(i)]{EGR}
For every nondegenerate collection $\{F_1, \ldots, F_n\}$ of elements of $A$, there exists a 
unique monic differential operator $L \in A[\partial]$ of order $n$ such that $L F_i =0$ for all $0 \leq i \leq n$. 
It is given by
\[
L g = W(F_1, \ldots, F_n, g)|_{n+1, n+1}, \quad \forall g \in A.
\]
\eth
The following theorem describes factorizations of such differential operators coming from
nondegenerate subcollections.
\bth{nondeg-fact} Assume that $\{F_1, \ldots, F_n\}$ is a nondegenerate collection of elements of $A$. 
Denote by $L \in A[\partial]$ the unique monic differential operator of order $n$ such that
\[
LF_i = 0, \forall 1 \leq i \leq n.
\]
Let $\{G_1, \ldots, G_k\}$ be a nondegenerate collection of elements of $A$ such that 
\[
L G_j = 0, \forall 1 \leq j \leq k.
\]
Denote by $P \in A[\partial]$ the unique monic differential operator of order $k$ such that 
\[
P G_i =0, \forall 1 \leq i \leq k.
\]
Then the following hold:
\begin{enumerate}
\item[(i)] \cite[Theorem 5.1]{K}
$L = Q P$ for some monic differential operator $Q \in A[\partial]$ of order $n-k$. 
\item[(ii)] If 
\[
G_1 = F_1, \ldots, G_k = F_k,
\]
then the collection $\{ P F_{k+1}, \ldots, P F_n \}$ of elements of $A$ is nondegenerate and 
$L = Q P$ where $Q \in A[\partial]$ is the unique monic differential operator of order $n-k$ such that 
$Q (P F_j) = 0$ for all $k + 1 \leq j \leq n$. 
\end{enumerate}
\eth
In \cite{K} Kasman proved a more general theorem than 
\thref{nondeg-fact}(i), which deals with operators with coefficients in an algebra $A$ that are formed by powers of a linear endomorphism of $A$ given by a sum of left multiplications and skew derivations. The special case of factorizations into linear terms was obtained earlier by Etingof, Gelfand and Retakh \cite[Theorem 1.1(b)]{EGR}; it has an additional condition that all subcollections $\{F_1, \ldots, F_l\}$ are nondegenerate for $1 \leq l \leq n$. We provide a proof of \thref{nondeg-fact}(i) for completeness. 

\begin{proof} (i) Since the operator $P \in A[\partial]$ is monic we can divide $L$ by it to get
\[
L= QP + R
\]
with a remainder 
\[
R = b_0 + \ldots + b_{k-1} \partial^{k-1} 
\]
with $b_0, \ldots b_{k-1} \in A$. Then $R G_j = 0$ for all $ 1 \leq j \leq k$, which is written in matrix form as
\[
(b_0, \ldots, b_{k-1}) W(G_1, \ldots, G_k) =0.
\]
Since the collection $\{G_1, \ldots, G_k\}$ is nondegenerate, $b_0 = \ldots = b_{k-1} = 0$. Therefore, $R=0$, and thus, $L = QP$. 

(ii) Denote the matrix
\[
B:=  \left(
\begin{matrix}
F_{k+1} & \ldots & F_n \\ 
\ldots & \ldots & \ldots \\ 
\delta^{k-1} F_{k+1} & \ldots &\delta^{k-1} F_n
\end{matrix} \right)
\in M_{k, n-k}(A).
\]
Consider the differential operator 
\[
\partial^j P = \partial^{j+k} + \sum_{i=0}^{j+k-1} a_{ij} \partial_i \in A [\partial]
\]
and the corresponding lower triangular matrix
\[
U := (u_{ij})_{i,j=1}^n \in M_n(A), \quad \mbox{where} \quad
u_{ij} := 
\begin{cases}
1, & \mbox{if} \; \; j=i \\
0, & \mbox{if} \; \; i < j \; \; \mbox{or} \; \; k \geq i > j \\
a_{ij}, & \mbox{if} \; \; i > j \; \; \mbox{and} \; \; i > k.
\end{cases}
\]
Since, $P F_i =0$ for all $1 \leq i \leq k$, we have
\begin{align}
\left(
\begin{matrix}
W(F_1, \ldots, F_k) & B \\
0 &W(PF_{k+1}, \ldots, P F_n) 
\end{matrix}
\right)
&= 
 \left(
\begin{matrix}
F_1 & \ldots & F_n \\ 
\ldots & \ldots & \ldots \\ 
\partial^{k-1} F_1 & \ldots & \partial^{k-1} F_n  \\ 
P F_1 & \ldots & P F_n  \\ 
\ldots & \ldots & \ldots \\ 
\partial^{n-k} P F_1 & \ldots &\partial^{n-k} P F_n
\end{matrix} \right)
\label{3matr}
\\
&= U W(F_1, \ldots, F_n),
\nn
\end{align}
where in the first term, $0$ denotes the zero matrix of size $(n-k) \times k$. 
The matrix $U \in M_n(A)$ is invertible because it is lower triangular matrix 
with 1's on the diagonal. Hence, the first matrix in \eqref{3matr} is invertible. 
That matrix has a block form, so its diagonal (block) entries must be invertible. 
Therefore, $W(PF_{k+1}, \ldots, P F_n) \in M_{n-k}(A)$ is invertible.  

Next we proceed with proving that $L = QP$ for the monic differential operator $Q \in A[\partial]$ defined in part (ii) 
of the proposition. 
Both differential operators $L$ and $QP$ are monic of order $n$ and annihilate the set $\{F_1, \ldots, F_n\}$. 
By \thref{EGR}, $L = QP$. 
\end{proof}
\subsection{Local picture for kernels of differential operators}
\label{3.2}
Let $\Om \subset \Cset$ be a domain (connected, simply connected open subset of $\Cset$). Denote 
\begin{itemize}
\item by $\MM_\Om(R)$ the algebra of $R$-valued meromorphic functions on $\Om$ and 
\item by $\DD_\Om(R)$ the algebra of $R$-valued meromorphic differential operators on $\Om$. 
\end{itemize}
We have the canonical embedding $R \hra \MM_\Om(R)$, given by constant functions,
which makes  $\MM_\Om(R)$ an $(\DD_\Om(R),R)$-bimodule. For $L(x, \partial_x) \in \DD_\Om(R)$, denote
\[
\ker L(x, \partial_x) := \{ f(x) \in \MM_\Om(R) \mid L(x, \partial_x) f(x) = 0 \}.
\]
Each $\ker L(x, \partial_x)$ has a canonical structure of right $R$-module. For a subdomain $\Om' \subseteq \Om$, set
\[
\ker_{\Om'} L(x, \partial_x) := \{ f(x) \in \MM_{\Om'}(R) \mid L(x, \partial_x) f(x) = 0 \}.
\]
The existence and uniqueness theorem for the initial value problem for vector-valued 
ordinary differential equations has the following corollary for the kernels of $R$-valued meromorphic differential operators:
\ble{ker-loc} Let $\Om$ be a domain in $\Cset$, $L \in \DD_\Om(R)$ and $x_0 \in \Om$ be a point at which $L$ is analytic. Then $x_0$ 
has a neighborhood $\OO_{x_0}$ such that $\ker_{\OO_{x_0}} L$ is a free {\em{(}}right\,{\em{)}} $R$-module of rank $n := \ord L$. 
More precisely, we have an isomorphism
of right $R$-modules
\[
\ker_{\OO_{x_0}} L \stackrel{\cong}\ra R^n, 
\] 
given by $\psi \in \ker_{\OO_{x_0}} L \mt (\psi(x_0), \psi'(x_0), \ldots, \psi^{(n-1)}(x_0))$.
\ele
\bde{nondeg} A free right $R$-submodule of $\MM_\Om(R)$ of finite rank will be called {\em{nondegenerate}} 
if one of its $R$-bases is nondenegerate in the sense of \deref{nondeg0}.
\ede
\bre{bases}
In the setting of the definition, all $R$-bases of the $R$-submodule of $\MM_\Om(R)$ in question 
are nondegenerate because, for each two bases, 
the Wronski matrix of one of them equals the product of the Wronski matrix of the other and an invertible matrix in $M_n(R)$.  
\ere

We will need the following easy lemma:
\ble{nondeg}
For a domain $\Om \subset \Cset$ and $F_1, \ldots, F_n \in \MM_\Om(R)$, the collection $\{F_1, \ldots, F_n\}$ is nondegenerate if and only if 
the Wronski matrix $W(F_1, \ldots, F_n)(x_0) \in M_n(R)$ 
is invertible for some $x_0 \in \Om$ at which $F_1, \ldots, F_n$ are holomorphic. 
\ele
\begin{proof} The left-to-right direction is obvious. For the converse statement, assume that $W(F_1, \ldots, F_n)(x_0) \in M_n(R)$ is invertible at some point 
$x_0 \in \Om$ at which $F_1, \ldots, F_n$ are holomorphic. Denote by $D$ the subset at which 
at least one of the functions $F_1, \ldots, F_n$ is not holomorphic. It is a discrete subset of $\Om$ because $\Om \subset \Cset$ 
is open and connected. For all points $x \in \Om \backslash D$, define the matrix 
\[
m(x) := W(F_1, \ldots, F_n)(x) \in M_n(R)
\]
and the determinant 
\[
d(x) := \det[ m(x) \cdot ]
\]
of the $\Cset$-linear endomorphism of $M_n(R)$ given by left multiplication by $m(x)$. The function $d(x)$ 
is holomorphic on $\Om \backslash D$ and is canonically extended to a $\Cset$-valued meromorphic function on $\Om$, 
which is not identically equal to zero because $d(x_0) \neq 0$ by \coref{regR} applied to the finite dimensional 
algebra $M_n(R)$. Denote by $D_0 \subset \Om \backslash D$ the zero set of 
$d(x)$. Since $\Om \backslash D$ is open and connected, $D_0$ is a discrete subset of $\Om \backslash D$. Applying \coref{regR} 
one more time, we obtain that for all $x \in \Om \backslash (D \cup D_0)$, the matrices $m(x) \in M_n(R)$ 
are invertible. Applying Kramer's rule for inverses of square matrices, gives that $m(x)^{-1}$ extends to an $M_n(R)$-valued meromorphic function on $\Om$.
\end{proof}
\subsection{Differential operators and kernels}
\label{3.3}
By the uniqueness of solutions of the initial value problem for vector-valued ordinary differential equations, for $L(x, \partial_x) \in \DD_\Om(R)$ and a
point $x_0 \in \Om$ at which $L(x, \partial_x)$ is analytic, the map
\begin{equation}
\label{eval}
\ker L \to R^n \quad \mbox{given by} \quad
F(x) \in \ker L \mt (F(x_0), F'(x_0), \ldots, F^{(k-1)}(x_0))
\end{equation}
is an injective homomorphism of right $R$-modules. In particular, 
\[
\dim_{\Cset} \ker L \leq \ord L \dim_{\Cset} R.
\]
When $L$ is restricted to a sufficiently small neighborhood of $x_0$, this becomes an isomorphism by \leref{ker-loc}. 

The next proposition describes the reconstruction of $R$-valued meromorphic differential operators from their kernels.
It is a refinement of \thref{EGR}. 

\bpr{diff-op} Let $\Om \subseteq \Cset$ be a domain and $k$ be a positive integer.
\begin{enumerate}
\item[(i)] Assume that $L(x, \partial_x) \in \DD_\Om(R)$, $\ord L =n$ and $\dim_{\Cset} \ker L= n \dim_{\Cset} R$ (the last condition 
can be always achieved by passing to a sufficiently small neighborhood of a point at which $L$ is analytic).
Then $\ker L$ is a nondegenerate free right $R$-module of rank $k$. 
\item[(ii)] Assume that $V \subset \MM_\Om(R)$ is a nondegenerate free right $R$-submodule of rank $n$. 
Then there exists a unique operator $L(x, \partial_x) \in \DD_\Om(R)$ such that $\ker L = V$. It is given by 
\begin{equation}
\label{QQ}
L(x, \partial_x) g(x) = |W(F_1, \ldots, F_n, g)|_{n+1, n+1} \quad \forall g(x) \in \MM_\Om(R),
\end{equation}
\end{enumerate}
where $F_1, \ldots, F_n$ is an $R$-basis of $V$.
\epr
\begin{proof} Part (i) follows by combining Lemmas \ref{lker-loc} and \ref{lnondeg}. 

(ii) By \reref{bases}, $\{F_1, \ldots, F_n\}$ is a nondegenerate collection of elements of $\MM_\Om(R)$. 
\thref{EGR} implies that the monic differential $L \in \DD_\Om(R)$ from \eqref{QQ} satisfies
\[
L(x, \partial_x) F_i(x) = 0, \quad \forall \; 1 \leq i \leq k.
\]
We have that $\ker L  \supseteq V$ and $\dim_{\Cset} V = n \dim_{\Cset} R 
\geq \dim_{\Cset} \ker L$. Therefore, $\ker L = V$. 

The uniqueness statement in part (ii) follows from \thref{EGR} because $\ker L = V$ implies $\{F_1, \ldots, F_n\} \in \Ker L$. 
\end{proof}

We note that the corresponding statement of \prref{diff-op} to the case of $R$-valued $C^m$ of $C^\infty$ differential operators 
on an open interval $I \subseteq \Rset$ also holds. 
For the proof of part (i) one cannot apply \leref{nondeg}. Instead one shows that  
\[
\det [W(F_1, \ldots, F_n)(x) \cdot] = - \det (a_{n-1}(x)) \det[ W(F_1, \ldots, F_n)(x) \cdot]
\]
in the notation \eqref{determinants}, 
where $L(x, \partial_x) = \partial_x^n + \sum_{i=0}^{n-1} a_i(x) \partial_x^i$, and from this one deduces that 
\[
\det[ W(F_1, \ldots, F_n)(x) \cdot] = \exp \big(- \int_{x_0}^x \det(a_{n-1}(t)) d t \big) \det[ W(F_1, \ldots, F_n)(x_0) \cdot].
\] 
Applying \coref{regR}, we obtain that, if $W(F_1, \ldots, F_n)(x)$ is invertible for one point in $I$, then it is invertible everywhere on $I$. 
\subsection{Factorizations of differential operators and kernels}
\label{3.4}
\bpr{fact} Assume that $\Om \subseteq \Cset$ is a domain and $L(x, \partial_x) \in \DD_\Om(R)$ is a monic differential operator 
of order $n$ such that 
\[
\dim_{\Cset} \ker L = n \dim_{\Cset} R
\]
(which can be always achieved by passing to a sufficiently small neighborhood of a point at which $L$ is analytic) 
and $1 \leq k \leq n$. 
\begin{enumerate}
\item[(i)] All factorizations 
\begin{equation}
\label{facto}
L(x, \partial_x) = Q(x, \partial_x) P(x, \partial_x) 
\end{equation}
into monic differential operators $Q, P \in \DD_\Om(R)$ of orders $n-k$ and $k$ 
are in bijection with the nondegenerate free right $R$-submodules of $\ker L$ of rank $k$.
This bijection is given by $P \mt \ker P$. The inverse bijection sends a free $R$-submodule $V$ of $\ker L$ to the 
unique monic differential operator $P$ with kernel $V$ given by \prref{diff-op}{\em{(}}ii{\em{)}}.
\item[(ii)] All nondegenerate free right $R$-submodules of $\ker L$ are direct summands.
\end{enumerate}
\epr  
We note that Goncharenko and Vesolov stated a factorization/intertwining theorem for monic differential operators in the special case of matrix valued differential functions, see \cite[Therem 1]{GV}. However, Theorem 1 in \cite{GV} does not require kernels to be nondegenerate; omitting this condition is not sufficient for forming factorizations.
\begin{proof} By \prref{diff-op} (i), $\ker L$ is a nondegenerate free right $R$-module of rank $n$. 

(I) Assume that we have a factorization of the form in \eqref{facto}. 
We will show that  $\ker P$ is a nondegenerate free $R$-submodule of $\ker L$ of rank $k$
which is a direct summand. 

Let $x_0 \in \Om$ be such that $P$ is analytic at $x_0$. \leref{ker-loc} implies that 
there is a neighborhood $\OO_{x_0}$ and $F_1(x), \ldots, F_k(x) \in \ker L$ such that 
\[
\ker_{\OO_{x_0}} P \; \; \mbox{is a free right $R$-module with basis} \; \; F_1(x)|_{\OO_{x_0}}, \ldots, F_k(x)|_{\OO_{x_0}}.
\]
It follows from \prref{diff-op} that $F_1(x)|_{\OO_{x_0}}, \ldots, F_k(x)|_{\OO_{x_0}}$ generate 
a nondegenerate free right $R$-submodule of $\ker_{\OO_{x_0}} L$, and by \leref{nondeg}, $F_1(x), \ldots, F_k(x)$ generate
a nondegenerate free right $R$-submodule of $\ker L$. Let $T \in \DD_\Om(R)$ be given by
\[
T(x, \partial_x) g(x) := |W(F_1, \ldots, F_k, g)|_{k+1, k+1} \quad \forall g(x) \in \MM_\Om(R). 
\]
It follows from \prref{diff-op} (ii) that the restrictions of $P$ and $T$ to $\OO_{x_0}$ are equal. Since the 
two operators have meromorphic coefficients and $\Om$ is open and connected, $P = T$. This means 
that $\ker P$ equals the nondegenerate free $R$-submodule of $\ker L$ with 
$R$-basis $F_1(x), \ldots, F_k(x)$. 

Consider the homomorphism of right $R$-modules
\[
\eta : \ker L /\ker P \to \ker Q \quad \mbox{defined by} \; \; g(x) \mt P(x, \partial_x) g(x). 
\]
This map is obviously injective and 
\[
\dim_{\Cset} (\ker L/\ker P) = (n-k) \dim_\Cset R, \quad \dim_\Cset Q \leq \ord Q \dim_{\Cset} R = (n-k) \dim_\Cset R.
\] 
Therefore, $\eta$ is an isomorphism and $\ker Q$ is a free right $R$-module (applying again \prref{diff-op} (i)). 
In particular, $\ker Q$ is projective and 
\[
\ker L \cong \ker Q \bigoplus \ker P
\]
as right $R$-modules.
This proves that $\ker P$ is a nondegenerate free $R$-submodule of $\ker L$ that is a direct summand.

(II) Now assume that $V$ is a nondegenerate free $R$-submodule of $\ker L$ of rank $k$.  
Choose bases
 $\{F_1, \ldots, F_n\}$ and $\{G_1, \ldots, G_k\}$ of $\ker L$ and $V$ respectively, 
 considered as right $R$-modules. Then $L(x, \partial_x)$ is the 
 unique monic differential operator in $\DD_\Om(R)$ of order $n$ such that
 $LF_i =0$ for all $1 \leq i \leq n$. 
 Let $P(x, \partial_x) \in \DD_\Om(R)$ be the unique monic differential operator of order $k$ such that 
$P F_i =0$ for all $1 \leq i \leq k$. Then $\Ker P = V$. 

\thref{nondeg-fact}(i) implies that $P(x, \partial_x)$ right divides $L(x, \partial_x)$ in  $\DD_\Om(R)$, 
yielding a factorization of the form \eqref{facto}.
 
The combination of (I) and (II) implies the statements (i) and (ii) of the proposition. 
\end{proof}

It is shown in \exref{important ex} that the inverse implication in part (ii) of \prref{fact} is false.
\bre{dir-summ} The proof of \prref{fact} gives that the {\em{free nondegenerate $R$-submodules of $\ker L$ of rank $k$}}
are precisely {\em{the $R$-submodules of $\ker L$ that have the form}} 
\begin{align*}
V = \{ F_1(x) \cdot R + \cdots + F_k(x) \cdot R \mid \, &\mbox{for some $R$-basis $F_1, \ldots, F_n$ of $\ker L$ such that}
\\
&\mbox{$F_1, \ldots, F_k$ is a nondegenerate collection} \}.
\end{align*}
This is a more elementary formulation similar to the scalar case $R= \Cset$ where one deals with bases of subspaces of vector spaces.
However, we note that the nondegeneracy condition is till present and the representation theoretic approach will play a key role in the following treatment. 
\ere
\sectionnew{A description of the decorated adelic Grassmannian $\wt{\Gr}^{\Ad}(R)$ and the fibers of $\pi : \wt{\Gr}^{\Ad}(R) \tra \Gr^{\Ad}(R)$}
\label{class}
In this section we prove Theorem A and characterize the fibers of the projection from the 
decorated adelic Grassmannian of $R$ to its adelic Grassmannian.
Throughout the section, we fix a complete set of primitive orthogonal idempotents of
the finite dimensional algebra $R$: 
\[
e_1, \ldots, e_m.
\]
Denote the multiplicities
\begin{equation}
\label{m(i)}
m(i) := \{ 1 \leq j \leq m \mid e_j R \cong e_i R \} 
\end{equation}
where $\cong$ stands for isomorphism of (right) $R$-modules.
\subsection{A reformulation of the definition of $\wt{\Gr}^{\Ad}(R)$}
\label{4.1}
Recall the definition of the decorated adelic Grassmannian of $R$, \deref{adGr}.
The next lemma provides a simplification of the definition of $\wt{\Gr}^{\Ad}(R)$, 
which will play a key role in  the classification of the points of $\wt{\Gr}^{\Ad}(R)$.
\ble{GrAd-simpl} The decorated adelic Grassmannian $\wt{\Gr}^{\Ad}(R)$ of the finite dimensional algebra $R$ consists of all
meromorphic functions on $\Cset \times \Cset$ of the form
\[
\Phi(x,z) = P(x, \partial_x) \exp(xz),
\]
where $P(x,\partial_x)$ is a monic $R$-valued differential operator on $\Cset$ with rational coefficients which right-divides 
a scalar differential operator $q(\partial_x)$, $q(t) \in \Cset[t] \backslash \{ 0 \}$.
\ele
\begin{proof} A monic polynomial $p(t) \in R[t]$ is a regular element of $R[t]$ by a leading term argument. 
\leref{det} implies that for every monic polynomial $p (t) \in R[t]$, there exists a polynomial $f(t) \in R[t]$ such that 
$p(t)f(t) = f(t) p(t) = q(t)$ for some monic polynomial $q(t) \in \Cset[t]$; 
by a leading term argument, the polynomial $f(t)$ would have to be monic too.  

Therefore, every monic $R$-valued differential operator on $\Cset$ with rational coefficients 
which right-divides $p(\partial_x)$ for some monic polynomial $p (t) \in R[t]$, 
also right-divides $q(\partial_x)$ for some polynomial $q(t) \in \Cset[t]$. 
The lemma now follows from \deref{adGr}.
\end{proof}
\subsection{Meromorphic right divisors}
\label{4.2}
Recall from the introduction that 
\begin{itemize}
\item $\MM(R)$ denotes the algebra of $R$-valued meromorphic functions on $\Cset$. 
\end{itemize}
Specializing the notation from \S \ref{3.1}, denote 
\begin{itemize}
\item by $\DD(R)$ be the algebra of $R$-valued meromorphic differential operators on $\Cset$.
\end{itemize}
The space  $\MM(R)$ has a canonical structure of $(\DD(R),R)$-bimodule.

Before proving Theorem A, we solve the simpler classification problem for those right-divisors of (nonzero)
scalar differential operators that are monic and meromorphic but do not necessarily have rational coefficients.

\bpr{class0} There is a bijection between the elements of 
\begin{enumerate}
\item the set of monic $R$-valued meromorphic differential operators on $\Cset$ {\em{(}}i.e., $P(x, \partial_x) \in \DD(R)${\em{)}} which right-divide 
a scalar differential operator $q(\partial_x)$, $q(t) \in \Cset[t] \backslash \{ 0 \}$ and
\item the set of nondegenerate finite rank free right $R$-submodules $V$ of  
\[
\QP(R) = \bigoplus_{\alpha \in \Cset} R[x] \exp(\alpha x).
\]
\end{enumerate} 
The bijection is given by $P \mt \ker P$. The inverse bijection sends a right $R$-submodule $V$ of $\QP(R)$ to the 
unique monic differential operator $P \in \MD(R)$ with kernel $V$ and order equal to the $R$-rank of $V$.

All submodules in (2) are direct summands of $\QP(R)$.
\epr
\begin{proof} Denote the two sets by $\AA$ and $\B$, respectively.

Each $P(x, \partial_x) \in \AA$, right-divides a monic differential operator of the form
\[
L(x, \partial_x) := \prod_{\al \in S} ( \partial_x - \al)^{N+1}
\]
for some finite subset $S \subset \Cset$ and $N \in \Nset$. The kernel of the latter is 
\begin{equation}
\label{qp-R-mod}
\ker L(x, \partial_x) = \sum_{\al \in S} (R + \cdots + R x ^N) \exp(\al x),
\end{equation}
and, by \prref{fact},

(*) $\ker P$ is a nondegenerate finite rank free right $R$-submodule of $\ker L$ which is a direct summand. 

Therefore, $\ker P \in \B$. In the opposite direction, assume that $V \in \B$. Then there exist a finite subset $S \subset \Cset$ and $N \in \Nset$, such that 
$V$ is a direct summand of the right $R$-module in \eqref{qp-R-mod}. \prref{diff-op}(ii) implies that 
there exists a unique monic differential operator $P \in \MD(R)$ of order equal to the $R$-rank of $V$ and 
such that $\ker P = V$. It follows from \prref{fact} that $P(x, \partial_x)$ right divides $L(x, \partial_x)$. 
Thus, $P \in \MD(R)$.  

The last statement of the proposition follows from (*). 
\end{proof}
\subsection{Proof of Theorem A, part I}
\label{4.3}
Recall from \S \ref{1.2} that $\RD(R)$ denotes the algebra of $R$-valued differential operators on $\Cset$ with rational coefficients.
In light of \leref{GrAd-simpl}, Theorem A is equivalent to the following:
\bth{A-2} Assume that $R$ is a finite dimensional complex {\em{(}}unital{\em{)}} algebra and that $e_1, \ldots, e_m$ is a complete set 
of primitive idempotents of $R$. For a free right $R$-submodule $V$ of $\QP(R)$ of $\rank_R V =l$, 
the following are equivalent:
\begin{enumerate}
\item $V$ is the kernel of a monic differential operator $P(x, \partial_x) \in \RD(R)$ of order $l$ that 
right-divides $q(\partial_x)$ for some $q(t) \in \Cset[t] \backslash \{0 \}$.
\item $V$ is a nondegenerate submodule of $\QP(R)$ and 
\[
V = \bigoplus_{\al \in \Cset} \left ( V \cap R[x] \exp(\al x) \right).
\]
\item $V$ is a nondegenerate submodule of $\QP(R)$ and 
\[
V = \bigoplus_{k=1}^n p_k(x) \exp(\al_k x) e_{i_k} R
\]
for some $p_1(x), \ldots p_n(x) \in R[x]$, $\al_1, \ldots, \al_n \in \Cset$ {\em{(}}repetitions allowed{\em{)}}.
\item $V$ is a nondegenerate submodule of $\QP(R)$ and has an $R$-basis of the form
\[
p_1(x) \exp(\alpha_1 x) e_1 + \cdots + p_m(x) \exp(\al_m x) e_m \quad \mbox{for some} \quad p_k(x) \in R[x], \alpha_k \in \Cset, i \leq m,
\]
where the exponents $\al_1, \ldots, \al_l$ are not necessarily distinct and different basis elements can share the same exponents. 
\end{enumerate}
All submodules $V$ of $\QP(R)$ with any of the above properties (1)--(4) are direct summands. 
\eth
In this subsection we deduce the equivalence $(1 \Leftrightarrow 2)$ from \prref{class0}.
In the next subsection we prove the equivalence $(2 \Leftrightarrow 3 \Leftrightarrow 4)$.
The last statement of the theorem follows from \prref{class0}.
\medskip
\\
\noindent
{\em{Proof of the equivalence $(1 \Leftrightarrow 2)$ in \thref{A-2}.}} 

$(1 \Rightarrow 2)$ Since $P(x, \partial_x) \in \RD(R)$, 
\begin{equation}
\label{Pprop}
P(x, \partial_x)  \big( R[x] \exp(\al x) \big) \subseteq R[x] \exp(\al x). 
\end{equation}
Let $f_1 + \cdots + f_n \in \ker P$ and $f_j \in R[x] \exp(\al_j x)$ for some distinct  
$\al_1, \ldots, \al_n \in \Cset$. Eq. \eqref{Pprop} and the fact that the sum \eqref{QP} is direct
imply that $f_j \in \ker P$ for all $1 \leq j \leq n$.
Therefore,
\[
\ker P = \bigoplus_{\al \in \Cset} (\ker P \cap R[x] \exp( \al x) ).
\]

$(1 \Leftarrow 2)$
Let $V$ be a nondegenerate free right $R$-submodule of $\QP(R)$ of $R$-rank $l$ and 
\begin{equation}
\label{Vcond}
V = \bigoplus_{\al \in \Cset} \left ( V \cap R[x] \exp(\al x) \right).
\end{equation}
\prref{class0} implies that there exists a unique monic differential operator 
\[
P(x, \partial_x) = \partial_x^l + a_1(x) \partial_x^{l-1} + \cdots + a_l(x) \quad 
\mbox{with} \quad a_k(x) \in \MM(R)
\]
such that $P(x,\partial_x)$ right-divides $q(\partial_x)$ for some $q(t) \in \Cset[t] \backslash \{0 \}$. 
It remains to be shown that the meromorphic functions $a_k : \Cset \to R$ are rational. It follows from 
\eqref{Vcond} that $V$ is spanned by functions of the form $b(x) \exp(\al x)$ 
for $b(x) \in R[x]$, $\al \in \Cset$. Therefore, the condition $V = \ker P$ is equivalent to the system
\[
(\partial_x - \al_j)^l b_j(x) + a_1(x) (\partial_x - \al_j)^{l-1} b_j(x) + \cdots + a_l(x), \quad 1 \leq j \leq N
\]
for some $N \geq l$, where $b_1(x), \ldots, b_N(x) \in R[x]$ and $\al_1, \ldots, \al_N \in \Cset$. 
We view this as a linear system of equations on the meromorphic functions $a_1, \ldots, a_l : \Cset \to R$.
By \prref{class0} we know that the system has a unique solution. Since the coefficients are polynomial $R$-valued 
functions, the unique solution is given by rational functions $a_1(x), \ldots, a_l(x)$. 
\qed
\subsection{Proof of Theorem A, part II}
\label{4.4}
The equivalence $(2 \Leftrightarrow 3 \Leftrightarrow 4)$ in \thref{A-2} is a consequence of the following more general 
result, which constintues the last statement in Theorem A. 
\bpr{A-3} For a finite rank free right $R$-submodule $V$ of $\QP(R)$, the following are equivalent:
\begin{enumerate}
\item $V = \bigoplus_{\al \in \Cset} \left ( V \cap R[x] \exp(\al x) \right)$.
\item $V = \bigoplus_{k=1}^n p_k(x) \exp(\al_k x) e_{i_k} R$
for some $p_1(x), \ldots p_n(x) \in R[x]$, $\al_1, \ldots, \al_n \in \Cset$ {\em{(}}repetitions allowed{\em{)}}.
\item $V$ has an $R$-basis of the form
\[
p_1(x) \exp(\alpha_1 x) e_1 + \cdots + p_m(x) \exp(\al_m x) e_m \quad \mbox{for some} \quad p_k(x) \in R[x], \alpha_k \in \Cset,
\]
where the exponents $\al_1, \ldots, \al_l$ are not necessarily distinct and different basis elements can share the same exponents. 
\end{enumerate}
\epr
For the proof of \prref{A-3} we will need the following lemma. We will denote by $R_R$ the right regular representation of $R$.
\ble{proj} If $U \subset R_R^{\oplus l}$ is an indecomposable projective $R$-submodule which is isomorphic to $e_i R$ (as a right $R$-module) 
for some $i \leq m$ and $l \in \Zset_+$, then there exist $r_1, \ldots, r_l \in R$ such that 
\[
U = (r_1, \ldots, r_l) e_i R.
\]
\ele
\begin{proof} Denote the isomorphism of right $R$-modules $\mu : e_i  R \stackrel{\cong}{\lra} U$. Set 
\[
(r_1, \ldots, r_l):=  \theta(e_i).
\]
Since $\theta(e_i) e_j = \theta(e_i e_j) = \delta_{j,i} \theta(e_i)$, we have $r_k e_j = \delta_{j,i} r_k$. Therefore, 
\[
r_k = r_k (e_1 + \cdots + e_m) = r_k e_i
\]
and $U = \theta(e_i) R = (r_1, \ldots, r_l) e_i R$.
\end{proof}
\noindent
{\em{Proof of \prref{A-3}.}} $(1 \Rightarrow 2)$ Let $V = \oplus_{\al \in \Cset} \left ( V \cap R[x] \exp(\al x) \right)$.
The assumption that $V$ is free implies that 
$V \cap R[x] \exp(\al x)$ is a projective right $R$-module. Since $V$ has finite rank, we can write
\begin{equation}
\label{decomp}
V \cap R[x] \exp(\al x) = V_{\al,1} \bigoplus \cdots \bigoplus V_{\al,s}
\end{equation}
for some indecomposable project right submodules $V_{\al,1}, \ldots, V_{\al,s}$. 
Hence, $V_{\al,k} \cong e_{i_k} R$ for some $i_k \leq m$.  \leref{proj} implies that 
\[
V_{\al,k} = p_k(x) \exp(\al x) e_{i_k} R \quad \mbox{for some} \; \; p_k(x) \in R[x]. 
\]
Finally, the module $V$ is the direct sum of modules modules of the form $V_{\al,k}$ for the exponents $\al \in \Cset$ such that 
\[
V \cap R[x] \exp(\al x) \neq 0.
\]

$(2 \Rightarrow 3)$ Assume that $V$ is a free right $R$-submodule of $\QP(R)$ of rank $l$ such that 
\begin{equation}
\label{V}
V = \bigoplus_{k=1}^n  p_k(x) \exp(\al_k x) e_{i_k} R
\end{equation}
with $p_1(x), \ldots p_n(x) \in R[x]$ and $\al_1, \ldots, \al_n \in \Cset$. Recall the definition \eqref{m(i)} of the 
multiplicities $m(i)$. Since the category of 
finite dimensional $R$-modules is Krull--Schmidt, the number of summands in \eqref{V}, 
isomorphic to $e_i R$, equals $l m(i)$. Since $e_1, \ldots, e_m$ is a complete set of primitive orthogonal idempotents 
of $R$, this gives a basis of $V$ satisfying the condition (3) by taking 
into account that 
\begin{multline*}
q_1(x) \exp(\al_1 x) e_1 R \bigoplus \cdots \bigoplus q_m(x) \exp(\al_m x) e_m R =  \\
\big( q_1(x) \exp(\al_1 x) e_1 + \cdots + q_m(x) \exp(\al_m x) e_m \big) R
\end{multline*}
is a free $R$-module for all $q_1(x), \ldots q_m(x) \in R[x]$ and $\al_1, \ldots, \al_m \in \Cset$. The 
basis consists of the elements $q_1(x) \exp(\al_1 x) e_1 + \cdots + q_m(x) \exp(\al_m x) e_m$ obtained in this way from 
\eqref{V}.

The implication $(3 \Rightarrow 1)$ follows from the identity
\[
(p_1(x) \exp(\alpha_1 x) e_1 + \cdots + p_m(x) \exp(\al_m x) e_m) e_i = p_i(x) \exp(\al_i x) e_i \quad \forall i \leq m.
\]
\qed
\bex{important ex} Consider the algebra
\[
R = M_2(\Cset)
\]
and the differential operator $\partial^3$ (with values in $M_2(\Cset)$), whose (matrix) kernel is
\[
\Ker \partial^3 = \{ p(x) \in M_2(\Cset)[x] \mid \deg p(x) \leq 2 \}. 
\]
Its vector kernel (in the space of meromorphic functions with values in $\Cset^2$) is the 6-dimensional vector space over $\Cset$ with basis
\[
\begin{pmatrix}
1 \\
0
\end{pmatrix}, \; \; 
\begin{pmatrix}
0 \\
1
\end{pmatrix}, \; \; 
\begin{pmatrix}
x \\
0
\end{pmatrix}, \; \; 
\begin{pmatrix}
0 \\
x
\end{pmatrix}, \; \; 
\begin{pmatrix}
x^2 \\
0
\end{pmatrix}, \; \; 
\begin{pmatrix}
0 \\
x^2
\end{pmatrix}.
\]
We obtain bases of the matrix kernel $\ker \partial^3$
as a right $M_2(\Cset)$-module by forming 3 matrices with columns given by the 6 vectors. Consider the basis
\[
F_1(x) =
\begin{pmatrix}
1 & 0 \\
0 &x
\end{pmatrix}, \quad 
F_2(x) =
\begin{pmatrix}
1 & 0 \\
0 &x^2
\end{pmatrix}, \quad 
F_3(x) =
\begin{pmatrix}
x & x^2 \\
0 & 0
\end{pmatrix} \in M_2(\Cset)[x] \subseteq \QP(M_2(\Cset)). 
\]
It shows that the right $M_2(\Cset)$-submodule 
\[
V:=F_1(x) M_2(\Cset) \oplus F_2(x) M_2(\Cset)
\]
is a direct summand of $\ker \partial^3$. At the same time 
the collection $\{F_1(x), F_2(x)\}$ is degenerate because
the first rows of $F'_1(x)$ and $F'_2(x)$ vanish.

Therefore the right $M_2(\Cset)$-submodule
\[
V \subset \partial^3
\]
is free and is a direct summand but is degenerate.
We also note that there are other ways to group the above 6 vectors into bases of $\partial^3$ giving rise to right $M_2(\Cset)$-submodules of $\ker \partial^3$ with the stated 3 properties. 

This shows that in Theorem A the condition that $V$ be nondegenerate cannot be replaced with the condition that $V$ be a direct summand and that the inverse implication in part (ii) of \prref{fact} is false. 
\eex
\subsection{The fibers of the projection $\wt{\Gr}^{\Ad}(R) \tra \Gr^{\Ad}(R)$} 
\label{4.5}
Recall the definitions of the decorated adelic Grassmannian $\wt{\Gr}^{\Ad}(R)$
and the adelic Grassmannain $\Gr^{\Ad}(R)$ from \S \ref{1.2}.
In Theorems \ref{tcover} and \ref{tdeco-to-ad} 
below, we explicitly describe the fibers of the projection 
\[
\pi : \wt{\Gr}^{\Ad}(R) \tra \Gr^{\Ad}(R), \quad \Phi(x,z) \mapsto \Phi_{\norm}(x,z)
\]
(recall \eqref{normal}) in terms of the classification of the points of $\wt{\Gr}^{\Ad}(R)$ from Theorem A. More precisely, in 
\thref{deco-to-ad} it is proved that these fibers consist of points of $\wt{\Gr}^{\Ad}(R)$ having a common successor with respect 
to a certain partial order and in \thref{cover} this partial order is explicitly described.
 
\bde{cover} Define a {\em{partial order}} on $\wt{\Gr}^{\Ad}(R)$ by setting
\[
\Phi(x,z) = P(x, \partial_x) \exp(xz) \prec \Psi(x,z) = S(x, \partial_x) \exp(xz) \in \wt{\Gr}^{\Ad}(R)
\]
if there exist scalars $\gamma_{kj} \in \Cset$, $1 \leq k \leq m$, $1 \leq j \leq l$ for some positive integer $l$ 
such that
\begin{equation}
\label{RtoS}
S(x, \partial_x) = P(x, \partial_x)  \prod_j \big( \partial_x- \gamma_{1j} e_1 - \cdots -  \gamma_{mj} e_m \big).
\end{equation}
\ede

The {\em{immediate successors}} of $\Phi$ are the functions $\Psi$ in the definition obtained for $j=1$.

Note that if $\Phi(x,z) = P(x, \partial_x) \exp(xz)  \in \wt{\Gr}^{\Ad}(R)$ and $S(x,\partial_x)$ is given by \eqref{RtoS}, then 
$\Psi(x,z) = S(x, \partial_x) \exp(xz) \in \wt{\Gr}^{\Ad}(R)$. The orthogonality of $e_j$ implies that $\Psi(x, \partial_x)$ is also given by
\[
\Psi(x, \partial_x) = P(x, \partial_x) \Big( \prod_j(\partial_x- \ga_{1j}) e_1 + \cdots +  \prod_j(\partial_x- \ga_{mj}) e_m \Big) \exp(xz).
\]

Our first goal is to describe this partial order explicitly. 
Let $A$ be a linear operator on a finite dimensional complex vector space $U$. There is a unique decomposition $U = U_0 \oplus U_1$ 
such that $A(U_i) \subseteq U_i$ for $i=0,1$ and 
\[
\mbox{$A|_{U_0}$ is nilpotent and $A|_{U_1}$ is invertible}.
\]
Define the $\Cset$-linear operator 
\[
F_A : U[x] \to U[x], \quad 
F_A(u(x)):= 
\begin{cases}
\exp(- A x) \int_0^x \exp(A t) u(t) dt, & \mbox{for} \; \; u(x) \in U_0[x] \\
\sum_{k=0}^\infty (-1)^{k+1} A^{-k -1} u^{(k)}(x), & \mbox{for} \; \; u(x) \in U_1[x].
\end{cases}
\]
In the first case $F_A(u(x)) \in U_0[x]$ because $\exp(x A)$ is polynomial in $x$; in the second case $F_A(u(x)) \in U_1[x]$ because
$u^{(k)}(x)=0$ for large $k$.  
\ble{int} For all $u(x) \in U[x]$, $F_A(u(x))$ is a solution of the differential equation 
\[
(\partial_x + A) F_A(u(x)) = u(x).
\]
\ele
\begin{proof} For $u(x) \in U_0[x]$ the statement is obvious. For $u(x) \in U_1[x]$, 
\[
F_A(u(x)) = \sum_{k=0}^\infty (-1)^{k+1} (A^{-1} \partial_x)^n (A^{-1} u(x))
\]
is a solution of $(A^{-1} \partial_x + 1) F_A(u(x)) =A^{-1} u(x)$ by using geometric series.
\end{proof}
\bth{cover} Assume that under the classification of Theorem A, $\Phi(x,z) \in \wt{\Gr}^{\Ad}(R)$
corresponds to 
\[
V := \bigoplus_{k=1}^n p_k(x) \exp(\al_k x) e_{i_k} R \subset \QP(R)
\]
for some $p_1(x), \ldots, p_n(x) \in R[x]$, $\al_1, \ldots, \al_n \in \Cset$ (repetitions allowed), where $n = m \ord P$. 

The immediate successors of $\Phi(x,z)$ in $\wt{\Gr}^{\Ad}(R)$ are parametrized by $\Cset^m$ and the one for 
$(\ga_1, \ldots, \ga_m) \in \Cset^m$ corresponds (under Theorem A) to
\[
V_{(\ga_1, \ldots \ga_n)}:= \big( \bigoplus_{i=1}^m \exp(\ga_i x) e_i R \Big) \bigoplus \Big( \bigoplus_{k=1}^n F_{A_{\al_k}} (p_k(x) e_{i_k}) \exp(\al_k x)  R \Big)
\subset \QP(R),
\]
where $A_{\al} \in \End_{\Cset} (R)$ is the operator of left multiplication by $(\al - \ga_1) e_1 + \cdots + (\al- \ga_m) e_m$. 
\eth
\begin{proof} Let $\Phi(x,z) = P(x, \partial_x) \exp(xz)$ with $P \in \RD(R)$, so $V = \Ker P$. The immediate successors of $\Phi(x,z)$ are 
the functions $\Phi\spcheck(x,z) = P\spcheck(x, \partial_x)  \exp(xz) \in \wt{\Gr}^{\Ad}(R)$ for
\[
P\spcheck(x, \partial_x) := P(x , \partial_x) (\partial_x- \gamma_1 e_1 - \cdots -  \gamma_m e_m) 
\]
and $(\ga_1, \ldots \ga_m) \in \Cset^m$. The proposition follows from the fact that
\[
\Ker P\spcheck = V_{(\ga_1, \ldots \ga_n)}.
\]
This equality holds because by \leref{int},
\[
(\partial_x- \gamma_1 e_1 - \cdots -  \gamma_m e_m ) V_{(\ga_1, \ldots \ga_n)} = \Ker P
\]
and, obviously, 
\[
\Ker (\partial_x- \gamma_1 e_1 - \cdots -  \gamma_m e_m )  \subset V_{(\ga_1, \ldots \ga_n)}.
\] 
\end{proof}
\bth{deco-to-ad} The fiber of $\pi : \wt{\Gr}^{\Ad}(R) \tra \Gr^{\Ad}(R)$ containing  $\Phi(x,z) \in \wt{\Gr}^{\Ad}(R)$ consists of all 
$\Psi(x,z) \in \wt{\Gr}^{\Ad}(R)$ such that $\Phi(x,z)$ and $\Psi(x,z)$ have a common successor.
\eth
\begin{proof} We need to prove that for $\Phi(x,z), \Psi(x,z) \in \wt{\Gr}^{\Ad}(R)$
\[
\Phi_{\norm} (x,z) = \Psi_{\norm}(x,z) \; \; \Leftrightarrow \; \; 
\exists \Theta(x,z) \in \wt{\Gr}^{\Ad}(R) \; \; 
\mbox{such that} \; \; \Phi \preceq \Theta \; \mbox{and} \; \Psi \preceq \Theta.  
\]
Let $\Phi(x,z) = P(x, \partial_x) \exp(xz)$ with $P \in \RD(R)$.
By Theorem A, $\Ker P$ has a basis of the form 
\[
p_{1j}(x) \exp(\alpha_{1j} x) e_1 + \cdots + p_{mj}(x) \exp(\al_{mj} x) e_m \quad \mbox{for} \quad p_{kj}(x) \in R[x], \alpha_{kj} \in \Cset, 1 \leq j \leq l,
\]
where $l = \ord P$. Denote $\underline{\al}:= (\al_{11}, \ldots, \al_{m1}, \ldots, \al_{1l}, \ldots, \al_{ml}) \in \Cset^{ml}$. 

($\Leftarrow$) It is sufficient to show that, if $\Phi\spcheck(x,z)$ is an immediate successor of $\Phi(x,z)$, then 
$\Phi\spcheck_{\norm}(x,z) = \Phi_{\norm}(x,z)$. Assume that, in the setting of  \thref{cover}, $\Phi\spcheck$ is obtained from $\Phi$ via the 
$m$-tuple $\underline{\ga}:=(\ga_1, \ldots, \ga_m) \in \Cset^m$. Set
\begin{align}
\label{g-al}
g_{\underline{\al}}(z) &:= \prod_{j=1}^l(z- \al_{1j}) e_1 + \cdots +  \prod_{j=1}^l(z- \al_{mj}) e_m \in R[z]
\\
\label{g-ga}
g_{\underline{\ga}}(z)&:= (z-\ga_1) e_1 + \cdots +  (z - \ga_m) e_m \in R[z].
\end{align}
Since $e_1, \ldots, e_m$ are orthogonal idempotents, 
\[
g_{\underline{\al}}(z) g_{\underline{\ga}}(z)= \Big[ (z-\ga_1) \prod_j(z- \al_{1j}) e_1 + \cdots +  (z - \ga_m) \prod_j(z- \al_{mj}) e_m \Big].
\]
\thref{cover} implies that
\begin{align*}
\Phi\spcheck_{\norm}(x,z) &= \Big[  P(x, \partial_x) (\partial_x- \gamma_1 e_1 - \cdots -  \gamma_m e_m ) \exp(xz) \Big] 
g_{\underline{\ga}}(z)^{-1} g_{\underline{\al}}(z)^{-1}
\\
&= \Big[  P(x, \partial_x) \exp(xz) \Big] g_{\underline{\ga}}(z)
g_{\underline{\ga}}(z)^{-1} g_{\underline{\al}}(z)^{-1} = \Phi_{\norm}(x,z).
\end{align*}

($\Rightarrow$) Let $\Psi(x,z) = S(x, \partial_x) \exp(xz)$ with $S \in \RD(R)$.
By Theorem A, $\Ker S$ has a basis of the form 
\[
q_{1j}(x) \exp(\be_{1j} x) e_1 + \cdots + q_{mj}(x) \exp(\be_{mj} x) e_m \quad \mbox{for} \quad q_{kj}(x) \in R[x], \be_{kj} \in \Cset, 1 \leq j \leq n,
\]
where $n = \ord S$. Set $\underline{\be}:= (\be_{11}, \ldots, \be_{m1}, \be_{1n}, \ldots, \be_{mn}) \in \Cset^{mn}$ and
\begin{equation}
\label{g-be}
g_{\underline{\be}}(z):= \prod_{j=1}^n (z- \be_{1j}) e_1 + \cdots +  \prod_{j=1}^n (z- \be_{mj}) e_m.
\end{equation}
The assumption $\Phi_{\norm} (x,z) = \Psi_{\norm}(x,z)$ 
implies
\[
\Big[  P(x, \partial_x) \exp(xz) \Big] g_{\underline{\al}}(z)^{-1} = \Big[  S(x, \partial_x) \exp(xz) \Big] g_{\underline{\be}}(z)^{-1}.
\]
Since $g_{\underline{\al}}(z) g_{\underline{\be}}(z) = g_{\underline{\be}}(z) g_{\underline{\al}}(z)$, we have
\[
\Big[  P(x, \partial_x) \exp(xz) \Big] g_{\underline{\be}}(z) = \Big[  S(x, \partial_x) \exp(xz) \Big] g_{\underline{\al}}(z).
\]
Therefore,
\[
P(x, \partial_x) g_{\underline{\be}}(\partial_x) \exp(xz)  = S(x, \partial_x) g_{\underline{\al}}(\partial_x) \exp(xz).
\]
Denote this function by $\Theta(x,z)$. Using one more time that $e_1, \ldots, e_m$ are orthogonal idempotents, we obtain
\begin{align*}
\Theta(x,z) &= P(x, \partial_x) \prod_{j=1}^n (\partial_x -\be_{1j} e_1 - \cdots - \be_{mj} ) \exp(x,z)
\\
&= S(x, \partial_x) \prod_{j=1}^l (\partial_x -\al_{1j} e_1 - \cdots - \al_{mj} ) \exp(x,z).
\end{align*}
This is an element of $\wt{\Gr}^{\Ad}(R)$ which is a common successor of $\Phi(x,z)$ and $\Psi(x,z)$.
\end{proof}
\sectionnew{A rational Grassmannian $\Gr^{\Rat}(R)$ and an embedding $\iota : \Gr^{\Ad}(R) \hra \Gr^{\Rat}(R)$}
\label{rational}
In this section we define a canonical embedding of the adelic Grassmannian $\Gr^{\Ad}(R)$ of a finite dimensional 
algebra $R$ in its rational Grassamannian $\Gr^{\Rat}(R)$ and characterize its image. In the special case 
$R = \Cset$, this recovers Wilson's results \cite{Wi1}. 
\subsection{A rational Grassmannian of $R[z]$}
\label{5.1} Recall from \S \ref{2.4} that the algebra $R(z)$ of $R$-valued rational functions in $z$ is also given by 
the localization
\[
R(z) = R[z] \big[ (\Cset[z] \backslash \{0 \})^{-1} \big].
\]

\bde{ratGr} Define {\em{the rational Grassmannian}} $\Gr^{\Rat}(R)$ 
to be the set of free left $R$-submodules $M \subset R(z)$ such that
\begin{enumerate}
\item[(i)] $h(z) R[z] \subseteq M \subseteq g(z)^{-1} R[z]$ for some $h(z), g(z) \in \Cset[z]\backslash \{0\}$ and
\item[(ii)]  $M/ \left( h(z) R[z] \right)$ is a direct summand of $g(z)^{-1} R[z] / \left( h(z) R[z]\right)$ as a left $R$-module
such that
\begin{equation}
\label{corank}
\rk_R \big( g^{-1}(z) R[z] / M \big) = \deg g(z).
\end{equation}
\end{enumerate}
\ede

\bre{freeratGr} \hfill
\begin{enumerate}
\item[(i)] 
\deref{ratGr} is equivalent to requiring that $M /(h(z) R[z])$ be a finite rank free left $R$-module instead of 
requiring the same condition on $M$. More specifically, in the presence of condition (i) and the first part of condition (ii), 
the condition that $M$ is free is equivalent to the condition that 
$M /(h(z) R[z])$ is a finite rank free left $R$-module. Indeed, since $g^{-1}(z) R[z] / \left( h(z) R[z]\right)$  is free, 
the first part of condition (ii) implies that $M/ \left( h(z) R[z] \right)$ is a projective $R$-module. Thus, 
\[
M \cong  h(z) R[z] \bigoplus M/ ( h(z) R[z] ).
\]
Therefore, $M$ is a free $R$-module if and only if $M/ \left( h(z) R[z] \right)$ is a free $R$-module, which is of finite rank due to condition (i). 

\item[(ii)] The first part of condition (ii) is equivalent to saying that $M$ is a direct summand of $g(z)^{-1} R[z]$.  

\item[(iii)]Equation \eqref{corank} is a corank condition on $M$ in the following sense. By the first part of condition (ii)
\[
g^{-1}(z) R[z] / \left( h(z) R[z]\right) \cong M/ \left( h(z) R[z] \right) \oplus N
\]
for a left $R$-module $N$. Since the left hand side is a finite rank free $R$-module and $M/ \left( h(z) R[z] \right)$
is a free $R$-module, so is $N$ by the Krull--Schmidt property. Equation \eqref{corank} requires that
\[
\rank_R N = \deg g(z). 
\]
\end{enumerate}
\ere

Recall the bilinear pairing 
\[
\lcor .,. \rcor :  R[z] \times \QP(R) \to R
\]
given by \eqref{pairing-formula}, and its properties \eqref{linearity}-\eqref{mixed-invar}. It follows from the second property that it descends to an $(R,R)$-bilinear map on the tensor product
\[
R[z] \otimes_R \QP(R) \to R.
\]
A stronger form of the invariance property \eqref{mixed-invar} is proved in \coref{invar} below.
\leref{pairing2} contains a second formulation of this pairing.

For a left $R$-module $W$, denote its dual module
\begin{equation}
\label{Dw}
D(W) := \Hom_R (W, R),  
\end{equation}
which is a right $R$-module via the action
\begin{equation}
\label{right-act}
(\nu .  r)(w) := \nu(w)r \quad \mbox{for} \quad \nu \in D(W), r \in R, w \in W. 
\end{equation}
It follows from \eqref{linearity} that the pairing \eqref{form} gives rise to the homomorphism of right $R$-modules
\[
\QP(R) \to D(R[z]), \quad f \mt \WW_f := \lcor - , f \rcor.
\]
In the scalar case when $R = \Cset$, the elements $\WW_f \in D(R[z])$ recover so called {\em{Wilson's conditions}} \cite{Wi1}.
For every subset $V \subset \QP(R)$, we have
\begin{equation}
    \label{Vperp}
V^\perp = \bigcap_{f \in V} \WW_f.
\end{equation}
\subsection{A characterization of the rational Grassmannian $\Gr^{\Rat}(R)$ in terms of submodules of $\QP(R)$}
\label{5.2}
\bth{bij} For every finite dimensional complex algebra $R$, there is a bijection between the following sets:
\begin{enumerate}
\item[(1)] the finite rank free right $R$-submodules $V$ of $\QP(R)$ which are direct summands
\end{enumerate}
and
\begin{enumerate}
\item[(2)] the free left $R$-submodules $N \subset R[z]$ which are direct summands and satisfy 
$h(z) R[z] \subseteq N$ for some $h(z) \in \Cset[z]$, $h(z) \neq 0$.
\end{enumerate}
The bijection is given by $V \subset \QP(R) \mt V^\perp \subseteq R[z]$. The reverse bijection is given by 
$M \subseteq R[z] \mt M^\perp \subset \QP(R)$. 
\eth
\bco{char-rat} Every element of the rational Grassmannian $\Gr^{\Rat}(R)$ has the form 
\begin{equation}
\label{perp-g}
g(z)^{-1} V^\perp
\end{equation}
for some finite rank free right $R$-submodule $V \subset \QP(R)$ which is a direct summand and $g(z) \in \Cset[z]$
such that $\deg g(z) = \rank_R V$. If $f_1, \ldots, f_l$ is an $R$-basis of $V$, then 
\[
V^\perp = \ker \WW_{f_1} \bigcap \ldots \bigcap \ker \WW_{f_l}.
\]

In the opposite direction, every element of the form \eqref{perp-g} belongs to $\Gr^{\Rat}(R)$. 
\eco
\begin{proof} If
\[
h(z) R[z] \subseteq M \subseteq g(z)^{-1} R[z]
\]
is an element of $\Gr^{\Rat}(R)$ as in \deref{ratGr}, then 
\[
N:= g(z) M \subset R[z]
\]
is a free left $R$-submodule which is a direct summand and satisfies $g(z) h(z) R[z] \subseteq N$. 
By \thref{bij}, there exists a finite rank free right $R$-submodule $V \subset \QP(R)$ which is a direct summand and $g(z) \in \Cset[z]$
such that 
\[
M= g(z)^{-1} V^\perp.
\]
From \leref{finitepairing} below it follows that 
\begin{equation}
\label{corank-perp}
\rk_R (R[z] /V^\perp) = \rk_R V.
\end{equation}
Thus, $\deg g(z) = \rk_R V$.

In the opposite direction, let $M:= g^{-1}(z) V^\perp$ be  of the form \eqref{perp-g}. Applying \thref{bij} in the opposite direction, 
gives that $N = g(z) M$ belongs to the set (2) in \thref{bij}. Therefore, $M$ is a free left $R$-module, satisfies
\[
h(z) R[z] \subseteq M \subseteq g(z)^{-1} R[z]
\]
for some $h(z) \in \Cset[z]$, and
$M/ \left( h(z) R[z] \right)$ is a direct summand of $g^{-1}(z) R[z] / \left( h(z) R[z]\right)$.
The corank condition \eqref{corank} follows from \eqref{corank-perp}.
\end{proof}

We start with some preparation for the proof of \thref{bij}.
Denote the complement to the diagonals in $\Cset^k$:
\[
\Cset^k \backslash \Delta := \{ (\al_1, \ldots \al_k) \in \Cset^k \mid \al_i \neq \al_j \; \; \mbox{for all} \; \; i \neq j \}.
\]
Fix $\vec{\al}:=(\al_1, \ldots, \al_k) \in \Cset^N \backslash \Delta$ and $\vec{N}:= (N_1, \ldots, N_k) \in \Nset^k$. Set
\begin{equation}
\label{h(z)}
h(z) := (z-\al_1)^{N_1} \ldots (z- \al_k)^{N_k} \in \Cset[z] .
\end{equation}
Consider the right $R$-submodule
\begin{equation}
\label{QPaN}
\QP_{\vec{\al}; \vec{N}}(R):= (R + \cdots + R x^{N_1-1}) \exp(\al_1 x) \bigoplus \cdots \bigoplus (R + \cdots + R x^{N_k-1}) \exp(\al_k x). 
\end{equation}
of $\QP(R)$.
\ble{finitepairing} The pairing \eqref{form} restricts to a perfect pairing between the finite rank free $R$-modules
\begin{equation}
\label{form2}
\lcor .,. \rcor :
R[z]/(h(z) R[z]) \times \QP_{\vec{\al}; \vec{N}}(R) \to R.
\end{equation}
\ele
As a restriction of the pairing \eqref{form}, the pairing in \leref{finitepairing} is left $R$-linear in the first argument 
and right $R$-linear in the second. We recall that a perfect pairing is one for which there exists a pair of bases $\{v_s\}$, $\{f_t\}$  
of the first (left) $R$-module and second (right) $R$-module which are dual with respect to the pairing:
\[
\lcor v_s, f_t \rcor = \delta_{st}.
\]
\begin{proof} It is clear that, with respect to the pairing \eqref{form}, 
\[
\QP_{\vec{\al}; \vec{N}}(R)^\perp = h(z) R[z].
\]
Thus, \eqref{form} restricts to a pairing between $\QP_{\vec{\al}; \vec{N}}(R)$ and 
$R[z]/(h(z) R[z])$. The tensor product decompositions
\[
R[z]/(q(z) R[z]) \cong R \otimes_\Cset \Cset[z]/(h(z) \Cset[z])
\quad \mbox{and} \quad
\QP_{\vec{\al}; \vec{N}}(R) \cong \QP_{\vec{\al}; \vec{N}}(\Cset) \otimes_\Cset R
\]
realize both modules as free left and right $R$-modules, respectively. A pair of dual bases of 
the two representations with respect to \eqref{form2} is obtained by using a pair of dual bases of 
the nondegenerate pairing of finite dimensional complex vector spaces
\[
\lcor .,. \rcor : \Cset[z]/(h(z) \Cset[z]) \times  \QP_{\vec{\al}; \vec{N}}(\Cset) \to \Cset
\]
where the pairing is given by \eqref{pairing-formula}.
\end{proof}
Denote by $(.)^\upVdash$ the orthogonal complement with respect to the pairing \eqref{form2}, while still denoting 
by $(.)^\perp$ the orthogonal complement with respect to the pairing \eqref{form}. Denote the projection 
\[
\kappa : R[z] \to R[z]/(q(z) R[z]). 
\]
\begin{proof}[Proof of \thref{bij}]
Let $V$ be an element of the set (2) in \thref{bij}. Then $V$ is a right $R$-submodule of $\QP_{\vec{\al}; \vec{N}}(R)$ for some 
$\vec{\al} \in \Cset^N \backslash \Delta$ and $\vec{N} \in \Nset^k$. Therefore 
\begin{equation}
V^\perp = \kappa^{-1} \big( V^\upVdash \big)
\label{ort1}
\end{equation}
and, in particular, $h(z) R[z] \subseteq V^\perp$, for $h(z) \in \Cset[z]$ given by \eqref{h(z)}.
Since $V$ is free and $V \oplus W = \QP_{\vec{\al}; \vec{N}}(R)$ for another right $R$-submodule $W$ of $\QP_{\vec{\al}; \vec{N}}(R)$,
$W$ is also free by the Krull--Schmidt property and 
\begin{equation}
\label{two-pairings}
V^\perp / (h(z) R[z]) \cong V^\upVdash \cong D(W)
\end{equation}
as left $R$-modules. In the first isomorphism we use \eqref{ort1} and in the second one the fact that the pairing \eqref{form2} is perfect.
Therefore, $V^\perp / (h(z) R[z])$ is a free left $R$-module of finite rank
(as the dual of a free right module of finite rank), and as a consequence, $V^\perp$ is also free. 
Moreover, \leref{finitepairing} implies that 
\[
V^\perp/(h(z) R[z]) \bigoplus W^\upVdash \cong R[z]/( h(z) R[z]).
\]
Hence, $V^\perp$ is an element of the set (2) in \thref{bij}.

Let $M$ be an element of the set (2) in \thref{bij}. So, 
$M \subset R[z]$ is a left $R$-submodule such that $M \supseteq h(z) R[z]$ 
for $h(z) \in \Cset[z]$. Assume that $h(z)$ is given by \eqref{h(z)}. Then
\begin{equation}
\label{ort2}
M^\perp = (M/(q(z) R[z]))^\upVdash.
\end{equation}
Since $M$ is free and $M/(q(z) R[z]) \oplus S = R[z]/( q(z) R[z])$ for a left $R$-submodule $S$ of $R[z]/(q(z) R[z])$, 
by the  Krull--Schmidt property, $S$ is free. Hence 
\[
M^\perp = (M/(q(z) R[z]))^\upVdash \cong D(S)
\]
is a free right $R$-submodule of $\QP_{\vec{\al}; \vec{N}}(R)$. (The last term refers to the dual of a left $R$-module.)
Finally, \leref{finitepairing} implies that
\[
M^\perp \bigoplus S^\upVdash = QP_{\vec{\al}; \vec{N}}(R).
\]
Hence, $M^\perp$ belongs to the set (1) in \thref{bij}.

It follows from \leref{finitepairing} and Eqs. \eqref{ort1}--\eqref{ort2} that in the settings of those equations
\[
\left(V^\perp \right)^\perp = V \quad \mbox{and} \quad
\left(M^\perp \right)^\perp = M,
\]
which completes the proof of \thref{bij}
\end{proof}
\subsection{An embedding of the adelic Grassmannian $\Gr^{\Ad}(R)$ into the  rational Grassmannian $\Gr^{\Rat}(R)$}
\label{5.3}
We first define a map
\[
\wt{\iota} : \wt{\Gr}^{\Ad}(R) \to \Gr^{\Rat}(R).
\]
Fix an element $\Phi(x,z) \in \Gr^{\Ad}(R)$, which under the classification of Theorem A corresponds to $V \subset \QP(R)$ 
with basis
\[
p_{1j}(x) \exp(\alpha_{1j} x) e_1 + \cdots + p_{mj}(x) \exp(\al_{mj} x) e_m \quad \mbox{for} \quad p_{kj}(x) \in R[x], \alpha_{kj} \in \Cset, 1 \leq j \leq l.
\]
Set
\begin{align}
\wt{\iota}(\Phi) &:= V^\perp . \Big[ \prod_j(z- \al_{1j}) e_1 + \cdots +  \prod_j(z- \al_{mj}) e_m \Big]^{-1} 
\label{iotatild}
\\
&= V^\perp . \Big[ \frac{\theta(z)}{\prod_j(z- \al_{1j})} e_1 + \cdots +  \frac{\theta(z)}{\prod_j(z- \al_{mj})} e_m \Big] . \theta(z)^{-1},
\nn
\end{align}
where
\[
\theta(z) = \prod_{ij} (z - \al_{ij}).
\]
The last expression in \eqref{iotatild} belongs to $\Gr^{\Rat}(R)$ because of \coref{char-rat} and \prref{pr-to-emb}, proved below.

\bth{ad-to-rat} The map $\wt{\iota} : \wt{\Gr}^{\Ad}(R) \to \Gr^{\Rat}(R)$ descends under $\pi : \wt{\Gr}^{\Ad}(R) \tra \Gr^{\Ad}(R)$ to an embedding 
$\iota : \Gr^{\Ad}(R) \hra \Gr^{\Rat}(R)$.

The image of the embedding $\iota : \Gr^{\Ad}(R) \hra \Gr^{\Rat}(R)$ consists of all points of the rational Grassmannian of the form 
\[
V^\perp = \Big(  \ker \WW_{f_1} \bigcap \ldots \bigcap \ker \WW_{f_l} \Big) . \Big[ \prod_j(z- \al_{1j}) e_1 + \cdots +  \prod_j(z- \al_{mj}) e_m \Big]^{-1}
\]
with $f_j(x) = p_{1j}(x) \exp(\alpha_{1j} x) e_1 + \cdots + p_{mj}(x) \exp(\al_{mj} x) e_m$ for some $p_{kj}(x) \in R[x]$, 
$\alpha_{kj} \in \Cset, 1 \leq j \leq l$.
\eth
First we derive several auxiliary results needed for the proof of the theorem.
The next lemma provides a second characterization of the pairing \eqref{pairing-formula}.

\ble{pairing2} The pairing $\lcor .,. \rcor :  R[z] \times \QP(R) \to R$ is also given by 
\[
\lcor p(z), f(x) \rcor =   p(\partial_x) f(x) |_{x=0}.
\]
\ele
\begin{proof}
It is sufficient to prove the formula for $p(z)  = sz^m$ and $f(x) = x^n \exp(\al x) r$ where $s, r \in R$, $\al \in \Cset$ and $m,n$ are nonnegative 
integers. By Leibnitz's rule we have 
\[
s\partial_x^m (x^n \exp(\al x)r)  = s \sum_{j=0}^{m} \binom{m}{j} n(n-1)\ldots (n-j+1)x^{n- j} \al^{m-j} \exp(\al x)r.
\]
When evaluated at $x=0$ the above sum has at most one summand, the one with $j= n$, when $m \geq n$ and has none when $m< n$. Hence 
\begin{align*}
&p(\partial_x) f(x) |_{x=0}=  s\partial_x^m (x^n \exp(\al x)r)|_{x=0}  = \binom{m}{n} n!  \al^{m-n} sr  
\\
&=  m(m-1)\ldots(m-n+1) \al^{m-n}  
= p^{(n)}(\al)r = \lcor p(z), f(x) \rcor.
\end{align*}
\end{proof} 
\bco{invar}
The pairing $\lcor .,. \rcor :  R[z] \times \QP(R) \to R$ satisfies the following invariance property:
\[
\lcor p(z) q(z) , f(x) \rcor = \lcor p(z), q(\partial_x) f(x) \rcor 
\]
for all $p(z), q(z) \in R[z]$ and $f(x) \in \QP(R)$. 
\eco
\begin{proof} Applying \leref{pairing2}, we obtain
\[
\lcor p(z) q(z) , f(x) \rcor = p(\partial_x) q(\partial_x) f(x)|_{x=0} = 
\lcor p(z), q(\partial_x) f(x) \rcor.
\]
\end{proof}
\bpr{pr-to-emb} Assume that $\Psi(x,z) \in \wt{\Gr}^{\Ad}(R)$ is an immediate successor of $\Phi(x,z) \in \wt{\Gr}^{\Ad}(R)$; that is 
\[
\Phi(x,z) = P(x, \partial_x) \exp(xz), \;  
\Psi(x,z) = S(x, \partial_x) \exp(xz), \; \; \mbox{and} 
\; \; S(x, \partial_x) = P(x, \partial_x) g(\partial_x)
\]
where $g(z) = (z- \ga_1)e_1 + \cdots + (z - \ga_m) e_m$ for some $\ga_1, \ldots, \ga_m \in \Cset$. If $V:= \Ker P$ and $W := \Ker S$, then  
\[
V^\perp g(z) = W^\perp.
\]
\epr
\begin{proof} Let $p(z) \in V^\perp$. Applying \coref{invar} and using that $V = g(\partial_x) W$, we obtain
\[
\lcor p(z) g(z), W \rcor = \lcor p(z), g(\partial_x) W \rcor = \lcor p(z), V \rcor =0.
\]
Hence, $p(z) g(z) \in W^\perp$, and thus, 
\[
V^\perp g(z) \subseteq W^\perp.
\]
There exist $\vec{\al}:=(\al_1, \ldots, \al_k) \in \Cset^N \backslash \Delta$ and $\vec{N}:= (N_1, \ldots, N_k) \in \Nset^k$
such that 
\[
V, W \subseteq \QP_{\vec{\al}; \vec{N}}(R),
\]
recall \eqref{QPaN}. In particular, we have that $\be_1, \ldots, \be_m \in \{ \al_1, \ldots, \al_k \}$. 
We will use the notation for $h(z)$ from \eqref{h(z)}. 
Since $\QP_{\vec{\al}; \vec{N}}(R)^\perp = h(z) R[z]$, 
\[
V^\perp, W^\perp \supseteq h(z) R[z].
\]
Because $g(z) | (z- \be_1) \ldots (z - \be_m)$, by making $N_1, \ldots, N_k$ larger, without loss of generality we can assume that 
\[
V^\perp g(z) \supseteq h(z) R[z].
\]
Consider the diagram
\[
\begin{tikzcd}
V^\perp g(z) \arrow[hookrightarrow]{r} 
\arrow[hookrightarrow]{dr}
& V^\perp \\
h(z) R[z] 
\arrow[hookrightarrow]{r}
\arrow[hookrightarrow]{u}
& W^\perp
\end{tikzcd}
\] 
Using \eqref{two-pairings} and the perfect pairing from \leref{finitepairing}, we get
\begin{align*}
\dim_\Cset V^\perp / (h(z) R[z]) &= \dim_\Cset V^\upVdash = \dim_\Cset QP_{\vec{\al}; \vec{N}}(R) - \dim_\Cset V, 
\\
\dim_\Cset W^\perp / (h(z) R[z]) &= \dim_\Cset W^\upVdash = \dim_\Cset QP_{\vec{\al}; \vec{N}}(R) - \dim_\Cset W 
\\
&=  \dim_\Cset QP_{\vec{\al}; \vec{N}}(R) - \dim_\Cset V - \dim_\Cset R.
\end{align*} 
We have $\dim V^\perp / (V^\perp g(z)) =  \dim_\Cset R$ because $g(z) \in R[z]$ is a monic polynomial of degree 1. 
This implies that $\dim_\Cset W^\perp / (h(z) R[z]) = \dim_\Cset (V^\perp g(z))/(h(z) R[z])$, and the 
embedding $V^\perp g(z) \subseteq W^\perp$ gives that $V^\perp g(z) = W^\perp$.
\end{proof}
\begin{proof}[Proof of \thref{ad-to-rat}] It follows from \prref{pr-to-emb} that, if $\Psi(x,z) \in \wt{\Gr}^{\Ad}(R)$ 
is an immediate successor of $\Phi(x,z) \in \wt{\Gr}^{\Ad}(R)$, then
\[
\wt{\iota} \circ \pi( \Phi) = \wt{\iota} \circ \pi( \Psi).
\] 
Iterating this, we obtain that the same equality holds if $\Psi$ is a successor of $\Phi$. 
\thref{deco-to-ad} implies that, if $\pi( \Phi) = \pi( \Psi)$ for $\Psi, \Phi \in \wt{\Gr}^{\Ad}(R)$, then 
$\wt{\iota} \circ \pi( \Phi) = \wt{\iota} \circ \pi( \Psi)$.  Hence, $\wt{\iota}$ descends to a map 
$\iota : \Gr^{\Ad}(R) \to \Gr^{\Rat}(R)$.

Next, we show that the map $\iota$ is an embedding. Assume that
\begin{equation}
\label{assumePf}
\wt{\iota} \circ \pi( \Phi) = \wt{\iota} \circ \pi( \Psi) \quad \mbox{for} \quad
\Phi(x,z), \Psi(x,z) \in \wt{\Gr}^{\Ad}(R).
\end{equation}
Write, $\Phi(x,z) = P(x, \partial_x) \exp(xz)$ and $\Psi(x,z) = S(x, \partial_x) \exp(xz)$
for some $P, S \in \RD(R)$. Furthermore,
\[
\Phi_{\norm}(x,z) = \Phi(x,z) g_{\underline{\al}} (z)^{-1}, \quad 
\Psi_{\norm}(x,z) = \Psi(x,z) g_{\underline{\be}} (z)^{-1}
\]
for some $g_{\underline{\al}} (z), g_{\underline{\be}} (z) \in R[z]$ given by \eqref{g-al} and \eqref{g-be} with
$\al_{ij}, \be_{ij} \in \Cset$. Set $V:= \Ker P$, $W:= \Ker S$. The assumption \eqref{assumePf}
means that 
\begin{equation}
\label{assume2}
V^\perp  g_{\underline{\al}} (z)^{-1} = W^\perp g_{\underline{\be}} (z)^{-1}.
\end{equation}
Consider the successors
\begin{align*}
\Phi\spcheck (x,z) &= P\spcheck(x,\partial_x) \exp(xz), \quad P\spcheck(x,\partial_x): = P(x, \partial_x)  g_{\underline{\be}} (\partial_x), \\
\Psi\spcheck (x,z) &= S\spcheck(x,\partial_x) \exp(xz), \quad S\spcheck(x,\partial_x): = S(x, \partial_x)  g_{\underline{\al}} (\partial_x),
\end{align*}
of $\Phi(x,z)$ and $\Psi(x,z)$, respectively. Denote $V\spcheck := \Ker P\spcheck$ and $W\spcheck := \Ker S\spcheck$. It follows 
from \prref{pr-to-emb} that 
\[
(V\spcheck)^\perp = V^\perp g_{\underline{\be}}(z) \quad \mbox{and} \quad
(W\spcheck)^\perp = W^\perp g_{\underline{\al}}(z).
\]
Combining this with \eqref{assume2} and using the commutativity of $g_{\underline{\al}}(z)$ and $g_{\underline{\be}}(z)$ (because $\{e_j\}$
are orthogonal idempotents) gives $(V\spcheck)^\perp = (W\spcheck)^\perp$, and thus, $V\spcheck= W\spcheck$. Therefore, 
$P\spcheck(x,\partial_x), S\spcheck(x,\partial_x) \in \RD(R)$ are monic differential operators with the same kernel. So,  
$P\spcheck(x,\partial_x) = S\spcheck(x,\partial_x)$, which in turn implies that $\Phi\spcheck (x,z) = \Psi\spcheck (x,z) \in \wt{\Gr}^{\Ad}(R)$
is a common successor of $\Phi(x,z)$ and $\Psi(x,z)$.  Invoking \thref{deco-to-ad} one more time, gives that 
$\pi(\Phi) = \pi (\Psi)$. This proves that $\iota$ is an embedding.

The last statement of the theorem about the image of $\iota$ follows from the definition of $\wt{\iota}$.
\end{proof}
\sectionnew{Examples}
\label{ex}
\bex{ge1} [{\em{Dual numbers}}] Consider the algebra of dual numbers
\[
R:= \Cset[\ep]/(\ep^2),
\]
see \exref{ex-fin-dim-alg}(ii). Its primitive idempotent is simply the element 1. Denote the elements 
\[
F_1 := e^x( 1 + (1+ \ep) x^2), \quad F_2:= e^x ( \ep + x) \in \QP(R),
\]
and the free right $R$-submodule of $\QP(R)$
\[
V:= F_1 R \oplus F_2 R.
\]
This module is nondegenerate: since $R$ is a commutative finite dimensional algebra, to check that the Wronski matrix of $F_1$ and $F_2$ is nondegenerate, all we need to show is that its determinant is not a zero divisor of $\QP(R)$. The latter is
\[
\det [ W(F_1,F_2)] = \det
\begin{pmatrix}
 F_1 & F_2 
 \\
 F'_1 & F'_2
\end{pmatrix}
= e^x( 1 - 2 \ep x - (1 + \ep) x^2).
\]
The right $R$-module $V \subset \QP(R)$ satisfies the conditions of Theorem A and the exponents in condition (ii) in this case are all equal to $1$.
The unique monic differential operator $P(x, \partial_x) \in \RD(R)$ with kernel $V$ is given by 
\[
P(x, \partial_x) = \partial_x^2 - 2\left( 
\frac{x^2 + x -1}{x^2 -1} - \ep\frac{x^2 +x +1}{(x^2-1)^2} \right) \partial_x + 
\frac{x+1}{x-1} - 2 \ep \frac{x-2}{(x-1)(x^2-1)}
\]
and the corresponding point of the decorated adelic Grassmmannian of $R$ is 
\[
\Phi(x, z) = P(x, \partial_x) e^{xz} \in \wt{\Gr}^{\Ad}(R) 
\]
and its image \eqref{normal} in the adelic Grassmannian is
\[
\Phi_{\norm} = \Phi(x,z)/(z-1)^2. 
\]
Recall the pairing \eqref{form}. The image of $\Phi_{\norm}$ 
in the rational Grassmannian of $R$ is given by  
\begin{align*}
\iota(\Phi_{\norm}) &= \frac{1}{(z-1)^2} V^\perp 
\\
&= \frac{1}{(z-1)^2} \{ f(z) \in R[z] 
\mid f(1) + (1 + \ep) f''(1) =0, \ep f(1) + f'(1)=0 \} 
\\
&= \frac{1}{(z-1)^2} \Span_R \{ 1 - \ep (z-1) + \frac{1-\ep}{2}(z-1)^2, 
(z-1)^3, (z-1)^4, \ldots \},
\end{align*}
where the elements in the last equation form an $R$-basis of 
$\iota(\Phi_{\norm})$. To verify that $\iota(\Phi_{\norm}) \in \Gr^{\Rat}(R)$ (see \deref{ratGr}) we note that
\[
(z-1) R[z] \subset \iota(\Phi_{\norm}) \subset (z-1)^{-2} R[z],
\]
and $\iota(\Phi_{\norm})$ is a direct summand of 
$(z-1)^{-2} R[z]$ with direct complement given by 
the free $R$-module with basis 
\[
\{ (z-1)^{-2}, (z-1)^{-1} \}.
\]

We have the surjective $\Cset$-algebra homomorphism $R \to \Cset$ given by $\ep \mt 0$. Under it, the above point in $\Gr^{\Rat}(R)$ goes to the point 
\[
\frac{1}{(z-1)^2} \Span_\Cset \{ 1 + (z-1)^2/2, 
(z-1)^3, (z-1)^4, \ldots \}
\]
of the classical Wilson's adelic Grassmannian $\Gr^{\Ad} \subset \Gr^{\Rat}$ obtained from the conditions $f(1) + f''(1) =0$ and 
$f'(1)=0$ supported at 1. The corresponding (normalized) wave function of the KP hierarchy is 
\[
\frac{z^2(x^2-1) - 2(x^2+x-1) z + (x+1)^2}{(z-1)^2(x^2-1)} e^{xz}.  
\]
\eex

\bex{ge2} [{\em{The path algebra of the $A_2$ quiver}}] Next we consider the case when $R$ is the path algebra of the quiver $1 \longrightarrow 2$. Here and in the next example, the paths will be written left to right as in \cite[Ch. II]{ASS}. The algebra is isomorphic to the algebra of lower triangular matrices
\[
\begin{pmatrix}
a & 0 
\\
b & c 
\end{pmatrix}, \quad a, b, c \in \Cset,
\]
see \cite[Example 1.3(c)]{ASS}. The primitive idempotents in the second picture are the elementary matrices $E_{11}$ and $E_{22}$. 
Take 
\[
F(x) := 
\begin{pmatrix}
(x+1) e^x & 0
\\
0 & x^2 -1
\end{pmatrix} \in \QP(R)
\quad \mbox{and} \quad 
V:= F(x) R \in \QP(R).
\]
The right $R$ submodule $V \subset QP(R)$ is nondegenerate since the condition of \leref{nondeg} are satisfied. 
Clearly it satisfies the conditions in Theorem A, but note that the exponents in front of different primitive idempotents in condition (ii) of Theorem A are not different (0 and 1), which illustrates the remark after Theorem A. 
The unique monic differential operator $P(x, \partial_x) \in \RD(R)$ is given by 
\[
P(x, \partial_x) = \partial_x - \begin{pmatrix}
\frac{x+2}{x+1} & 0
\\
0 & \frac{2x}{x^2-1}
\end{pmatrix}.
\]
The corresponding point of the decorated adelic Grassmannian of $R$ is
\[
\Phi(x,z) = P(x, \partial_x) e^{xz} = 
z e^{xz} - 
\begin{pmatrix}
\frac{x+2}{x+1} & 0
\\
0 & \frac{2x}{x^2-1}
\end{pmatrix} e^{xz} \in \wt{\Gr}^{\Ad}(R)
\]
and the point of the adelic Grassmannian is
\[
\Phi_{\norm}(x,z) = \Phi 
\begin{pmatrix}
(z-1)^{-1} & 0
\\
0 & z^{-1}
\end{pmatrix} \in \Gr^{\Ad}(R).
\]
The orthogonal complement $V^\perp$ with respect to the pairing \eqref{form} consists of the matrices of the form 
\[
\begin{pmatrix}
p_{11}(z) & 0
\\
p_{21}(z) & p_{22}(z) 
\end{pmatrix}
\]
such that
\[
p'_{11} (1) + p_{11}(1) =
p'_{21}(1) + p_{21}(1) = p''_{22}(0) - p_{22}(0) =0.
\]
A direct calculation shows that 
$V^\perp \begin{pmatrix}
z & 0 
\\
0 & z-1
\end{pmatrix}$ consists of the matrices of the form
\[
\begin{pmatrix}
c_{11} z^2 + z(z-1)^2 q_{11}(z) & 0 
\\
c_{21} z^2 + z (z-1)^2 q_{21}(z) 
& 
c_{22} (z-1)(z^2 + 2)/2 
+ c'_{22} z (z-1) + z^3(z-1) q_{22(z)}
\end{pmatrix},
\]
where $c_{ij}, c'_{ij} \in \Cset$, $q_{ij}(z) \in \Cset[z]$.
From this one obtains that the image of 
$\Phi_{\norm} \in \wt{\Gr}^{\Ad}(R)$ under $\iota$ is
\[
\iota(\Phi_{\norm})=
V^\perp \big( (z-1) E_{11} + z E_{22} \big)^{-1}
=
\frac{1}{z(z-1)} V^\perp 
\begin{pmatrix}
z & 0 
\\
0 & z-1 
\end{pmatrix}
\subset z^{-1}(z-1)^{-1} R[z]
\]
and that this is a free left $R$-module with basis 
\[
\begin{pmatrix}
z/(z-1) & 0 
\\
0 & (z^2+2)/2z
\end{pmatrix}, \quad 
\begin{pmatrix}
z-1 & 0
\\
0 & 1
\end{pmatrix}
, \quad 
\begin{pmatrix}
z^k (z-1) & 0
\\
0 & z^{k+1}
\end{pmatrix}, 
\; k \geq 1.
\]
To verify that $\iota(\Phi_{\norm}) \in \Gr^{\Rat}(R)$ (see \deref{ratGr}), we note that
\[
z^2 (z-1) R[z] \subset
\iota(\Phi_{\norm}) \subseteq z^{-1} (z-1)^{-1} R[z]
\]
and $\iota(\Phi_{\norm})$ is a direct summand of 
$z^{-1} (z-1)^{-1} R[z]$ with direct complement given by the free left $R$-module of rank 
$2 = \deg z (z-1)$ with basis  
\[
\{ z^{-1}(z-1)^{-1},  z^{-1} \}.
\]
considered as scalar elements of $R[z]$.
\eex

\bex{ge3} [{\em{The path algebra of the second  Kronecker quiver}}] The Kronecker algebra is the path algebra of the quiver 
\[
\begin{tikzcd}
1 
\arrow[r,bend right] \arrow[r,bend left] 
& 2
\end{tikzcd}
\]
The algebra has a  matrix presentation in the form
\[
\begin{pmatrix}
a & 0 & 0\\
b& d &  0\\
c& 0 & d
\end{pmatrix}, 
\]
where $a,b, c, d \in \Cset$. The primitive idempotents in the second presentation are $E_{11}$ and $E_{22} + E_{33}$. Take 
\[
 F_1(x) = \begin{pmatrix}
1& 0 & 0\\
1&x^2   &  0\\
 0& 0 &x^2
 \end{pmatrix}, \quad F_2(x) = \begin{pmatrix}
 x^2& 0 & 0\\
0&x   &  0\\
 1& 0 &x
 \end{pmatrix}
\in R[x] \subset \QP(R) 
\]
and 
\[
V:= F_1(x) R \oplus F_2(x) R \subset R[x] \subset \QP(R).
\]
It is straightforward to check that $V$ is a free right $R$-module of rank 2 and is nondegenerate. Clearly it satisfies the conditions in Theorem A and all exponents in condition (ii) of Theorem A equal 0.
The corresponding point of the decorated adelic Grassmannian of $R$ is
\[
\Phi(x,z) = P(x, \partial_x) e^{xz}, 
\]
where

 \[
 P(x, \partial_x) = \partial_x^2 - (F_1'', F_2'') \begin{pmatrix}
 F_1& F_2\\
 F_1'& F_2'
 \end{pmatrix}^{-1} \begin{pmatrix}
1\\
\partial_x
 \end{pmatrix}.  
 \]
 After some linear algebra we obtain that
 \[
 \begin{pmatrix}
 F_1& F_2\\
 F_1'& F_2'
 \end{pmatrix}^{-1}  = \begin{pmatrix}
 \al& \be\\
 \ga&\de
 \end{pmatrix}
 \]
 with  
\begin{align*}
 \al &=  \begin{pmatrix}
 1&0& 0\\
\frac{-1}{x(x-2)}   &\frac{1}{x(x-2)}  & 0\\
 0&0  & \frac{1}{x(x-2)}  
 \end{pmatrix},  &\be &= \begin{pmatrix}
\frac{-x}{2}  &0& 0\\
 \frac{1}{2x}&   \frac{1}{x}  & 0\\
 0&0 &   \frac{1}{x}
 \end{pmatrix},\\
  \ga &= \begin{pmatrix}
 0&0& 0\\
 \frac{-2}{x-2}  &   \frac{2}{x-2}   & 0\\
0  &0 & \frac{2}{x-2} 
 \end{pmatrix}, 
 &\de &=  \begin{pmatrix}
  \frac{-1}{2x}  &0& 0\\
 0 & -1& 0\\
 0  &0 & -1
 \end{pmatrix}.
\end{align*}
This gives the explicit form of the operator $P(x, \partial_x)$:
\[
P(x,\partial_x) = \partial_x^2 -   \frac{1}{x} \begin{pmatrix}
-1& 0& 0\\
1 & 2& 0\\
0 & 0& 2
\end{pmatrix}
\partial_x   +  \frac{1}{x(x-2)}
\begin{pmatrix}
(x-2)  &0& 0\\
1   & -1  & 0\\
0&0  & -1  
\end{pmatrix}.
\]
The corresponding point \eqref{normal} of the adelic Grassmannian is 
\[
\Phi_{\norm}(x,z)
= e^{xz} -   \frac{e^{xz}}{xz} \begin{pmatrix}
-1& 0& 0\\
1 & 2& 0\\
0 & 0& 2
\end{pmatrix}
+  \frac{e^{xz}}{x(x-2)z^2}
\begin{pmatrix}
(x-2)  &0& 0\\
1   & -1  & 0\\
0&0  & -1  
\end{pmatrix}.
\]
\eex
By a direct calculation one shows that its image under 
$\iota$, $\iota(\Phi_{\norm}) = z^{-2} V^\perp$, is the left $R$-submodule of $R(z)$ consisting of functions of the form
\[
z^{-2}
\begin{pmatrix}
a(z) & 0 & 0\\
b(z)& d(z) &  0\\
c(z)& 0 & d(z)
\end{pmatrix}, 
\]
where $a(z), b(z), c(z), d(z) \in \Cset[z]$ and 
\[
a(0) = a''(0) =0, \; \; 
d'(0) = d''(0) = 0, \; \; 
c(0) = 0, \; \; b''(0) =0, \; \; 
b(0) = c''(0) = - d(0). 
\]
It is a free left $R$-module with basis 
\[
\begin{pmatrix}
 z^{-1}     & 0      &  0      \\
-z^{-2}     & z^{-2} &  0      \\
-1          & 0      & z^{-2}
\end{pmatrix}, 
\quad 
\begin{pmatrix}
z^k & 0       &  0     \\
0       & z^k &  0     \\
0       & 0       & z^k
\end{pmatrix}, 
\; \; k \geq 1.
\]
We note the interesting feature of this example, compared to the previous one, that there is no $R$-basis 
of $\iota(\Phi_{\norm})$ consisting of diagonal matrices. 
To verify that $\iota(\Phi_{\norm}) \in \Gr^{\Rat}(R)$ (see \deref{ratGr}), we check that
\[
z  R[z] \subset
\iota(\Phi_{\norm}) \subseteq z^{-2} R[z]
\]
and $\iota(\Phi_{\norm})$ is a direct summand of 
$z^{-2} R[z]$ with direct complement given by the free left $R$-module of rank 
$2 = \deg z^2$ with basis  
\[
\begin{pmatrix}
z^{-2} & 0 & 0 \\
0 & z^{-1} & 0 \\
0 & 0 & z^{-1}
\end{pmatrix}, 
\quad 
\begin{pmatrix}
1 & 0 & 0 \\
0 & 1 & 0 \\
0 & 0 & 1
\end{pmatrix}. 
\]

\end{document}